\documentclass[11pt]{article} 
\usepackage{amsfonts,amsmath,amssymb}
\usepackage{color}
\usepackage{amsthm}
\usepackage{colortbl}
\usepackage{enumerate}
\usepackage{a4wide}
\usepackage{graphicx}
\usepackage{url}
\usepackage{subfigure}
\usepackage{xcolor}
\usepackage{dsfont}
\usepackage{fancybox}
\usepackage[small,nohug,heads=littlevee]{diagrams} 
\diagramstyle[labelstyle=\scriptstyle]

\newcommand{\R}{\mathbb{R}}
\newcommand{\E}{\mathbb{E}}

\newcommand\al{\alpha}

\newcommand\pref[1]{(\ref{#1})}

\newtheorem{remark}{Remark}
\newcommand{\bfm}[1]{\boldsymbol{#1}}

\newcommand{\fa}{\mathfrak{a}}
\newcommand{\fc}{\mathfrak{c}}
\newcommand{\fm}{\mathfrak{m}}
\newcommand{\fd}{\mathfrak{d}}
\newcommand{\fb}{\mathfrak{b}}
\newcommand{\fq}{\mathfrak{q}}
\newcommand{\fr}{\mathfrak{r}}

\DeclareMathOperator{\argmin}{argmin}

\def\<#1,#2>{\left<#1,#2\right>}
\let\bar\overline

\newcommand{\bone}{{\bf 1}}

\newcommand{\EE}{\mathbb{E} }

\newcommand{\PP}{\mathbb{P} }

\newcommand{\RR}{\mathbb{R} }

\newcommand{\cH}{\mathcal{H} }

\newcommand{\cL}{\mathcal{L} }
\newcommand{\cP}{\mathcal{P}}

\def\cA{\mathcal A}

\def\t{\tilde}
\def\o{\overline}
\def\u{\underline}

\def\reff#1{{\rm(\ref{#1})}}

\newtheorem{theorem}{Theorem}[section]

\newtheorem{proposition}[theorem]{Proposition}

\title{ \textbf{Control of McKean-Vlasov Dynamics versus Mean Field Games}\footnote{Paper presented at the  conference "Stochastic Partial Differential Equations: Analysis, Numerics, Geometry and Modeling", ETH Z{\"{u}}rich, September 12-16, 2011  and the Third Humboldt-Princeton Conference, October 28-29, 2011}
}
\author{Ren\'e Carmona${}^{(a),}$\footnote{rcarmona@princeton.edu, Partially supported  by NSF: DMS-0806591 } \hspace*{.05cm},
Fran{\c c}ois Delarue${}^{(b),}$\footnote{delarue@unice.fr} \hspace*{.05cm} and
Aim\'e Lachapelle${}^{(a),}$\footnote{alachape@princeton.edu, Partially supported  by NSF: DMS-0806591}
 \\ \\
(a) ORFE, Bendheim Center for Finance, Princeton University,
\\
Princeton, NJ  08544, USA.
\\
\\
(b) Laboratoire J.A. Dieudonn\'e, Universit\'e de Nice Sophia-Antipolis, 
\\
Parc Valrose, 06108 Nice Cedex 02, France.}

\date{}
\begin{document}
\maketitle

\begin{abstract}
We discuss and compare two methods of investigations for the asymptotic regime of stochastic differential games with a finite number of players as the number of players tends to the infinity. These two methods differ in the order in which optimization and passage to the limit are performed. When optimizing first, the asymptotic problem is usually referred to as a mean-field game. Otherwise, it reads as an optimization problem over controlled dynamics of McKean-Vlasov type. Both problems lead to the analysis of forward-backward stochastic differential equations, the coefficients of which depend on the marginal distributions of the solutions. We explain the difference between the nature and solutions to the two approaches by investigating the corresponding forward-backward systems. General results are stated and specific examples are treated, especially when cost functionals are of linear-quadratic type. \end{abstract}
\section{Introduction}
\label{se:introduction}

The problem studied in this paper concerns stochastic differential games with a large number of players. We compare two methods of investigation which offer, in the asymptotic regime when the number of players tends to infinity, a structure which is simple enough to be amenable to actual solutions, both from the theoretical and  numerical points of view.

In order to derive tractable solutions, we assume that all the players are similar in their behavior, and that each individual on his own, can hardly influence the outcome of the game. We further strengthen the symmetry of the problem by assuming that the interaction between the players is of \emph{mean field type} in the sense that whenever an individual player has to make a decision, he or she sees only averages of functions of the private states of the other players. These games are symmetric in the sense that all the players are statistically identical, but they are not anonymous\footnote{In the anonymous framework, dynamics are given for the statistical distribution of the population so that the private dynamics of the players are not explicit.} (see for example \cite{JovanovicRosenthal}) or with weak interaction in the sense of \cite{Horst2}. In the large game limit, a given player should feel the presence of the other players through the statistical distribution of the private states of the other players, and should determine his optimal strategy by optimizing his appropriately modified objective criterion taking the limit $N\to\infty$ into account.
The search for an approximate Nash equilibrium of the game in this asymptotic regime is dictated by the Mean Field Game (MFG for short) proposition of Lasry and Lions. See for example \cite{MFG1,MFG2,MFG3,MFG4,Lions,GueantLasryLions.pplnm}.

Without the optimization part, and when all the players use the same distributed feedback strategy as suggested by the symmetry of the set-up, the large population regime is reminiscent of Mark Kac's propagation of chaos theory put on rigorous mathematical ground by McKean and known under the name of  McKean--Vlasov (MKV for short) theory. See for example Sznitman's beautiful mathematical treatment 
\cite{Sznitman}. Indeed, according to these works, one expects that, in the limit $N\to\infty$, the private states of the individual players evolve independently of each other, each of them satisfying a specific stochastic differential equation with coefficients depending upon the statistical distribution of the private state in question. Having each player optimize his own objective function under the constraints of this new private state dynamics amounts to the stochastic control of McKean--Vlasov dynamics, mathematical problem not understood in general. See nevertheless \cite{AndersonDjehiche} for an attempt in this direction. Coming back to our original problem we may wonder if the optimization over the feedback strategies in this new control problem leads to one form of approximate equilibrium for the original $N$ player game. In any case, if such a problem can be solved, the issue is to understand how an optimal feedback strategy for such a problem compares with the result of the MFG analysis.

The main thrust of this paper is to investigate the similarities and the differences between the MFG approach to stochastic games with mean field interactions, and the search for an optimal strategy for controlled McKean--Vlasov stochastic differential equations. In Section \ref{se:games}, we give a pedagogical introduction to these two very related (and sometimes confused) problems with the specific goal to explain the differences between their nature and solutions; in Section \ref{se:solvability}, we also provide the reader with a short guide for tackling these problems within the framework of Forward--Backward Stochastic Differential Equations (FBSDEs) of the McKean--Vlasov type. 

The remaining part of the paper is devoted to the discussion of special classes of models for which we can address the existence issue and compare the solutions of the two problems explicitly. We give special attention to Linear Quadratic (LQ) stochastic games because of their tractability. Indeed, in this case, the mean field character of the interaction is of a very simple nature as it involves only the empirical mean and the empirical variance of the individual states, and they both enter the coefficients linearly. For these models, we implement the MFG approach and analyze the optimal control of MKV dynamics in Sections \ref{se:LQMFG} and \ref{se:LQMKV}\footnote{After completion of this work, we were made aware of the appearance on the web of a very recent technical report  by A. Bensoussan, K. C. J. Sung,  S. C. P. Yam, and S. P. Yung entitled \emph{Linear Quadratic Mean Field Games}. In this independent work, the authors present a study of linear quadratic mean field games in relation to the control of McKean-Vlasov dynamics very much in the spirit of what we do in Section \ref{se:LQMFG} and Section \ref{se:LQMKV}.} respectively.  While bearing some similarity to the contents of this section, the results of \cite{Bardi} on some linear quadratic MFGs are different in the sense that they concern some infinite horizon stationary cases without any attempt to compare the results of the MFG approach to the control of the corresponding 
McKean--Vlasov dynamics. We provide more  explicit examples in Section \ref{se:examples}, including a simple example of CO${}_2$ emissions regulation which serves as motivation for models where the interaction is of mean field type and appears only in the terminal cost. 

While the MFG approach does not ask for the solution of stochastic equations of the McKean--Vlasov type at first, the required fixed point argument identifies one of these solutions of standard optimal control problems  as the \emph{de facto} solution of an FBSDE of the McKean--Vlasov type as the marginal distributions of the solution appear in the coefficients of the equation.
Since the McKean--Vlasov nature of these FBSDEs is rather trivial in the case of LQ mean field games, we devote Section \ref{se:solvability} to a discussion of the construction of solutions to these new classes of FBSDEs which, as far as we know, have not been studied in the literature.

\section{Stochastic Differential Game with Mean Field Interactions}
\label{se:games}
We consider a class of stochastic differential games where the interaction between the players is given in terms of functions 
of average characteristics of the private states and actions of the individual players, hence their name Mean Field Games.
In our formulation, the state of the system (which is controlled by the actions of the individual players) is given at each time $t$ by a vector $X_t=(X^1_t,\cdots,X^N_t)$ whose $N$ components $X^i_t$ can be interpreted as the \emph{private states} of the individual players. 
A typical example capturing the kind of symmetry which we would like to include is given by models in which the dynamics of the private states are given by coupled stochastic differential equations of the form
\begin{equation}
\label{fo:firstd}
dX^i_t=\frac1N\sum_{j=1}^N\tilde b(t,X^i_t,X^j_t,\alpha^i_t)dt + \sigma dW^i_t,
\end{equation}
where $\tilde b$ is a function of time, the values of two private states, and the control of one player, and $((W^i)_{t \geq 0})_{1 \leq i \leq N}$ are independent Wiener processes. 
For the sake of simplicity, we assume that each process $(X^i_t)_{0\le t\le T}$ is univariate. Otherwise, the notations become more involved while the results remain essentially the same. The present discussion can accommodate models where the volatility $\sigma$ is a function with the same structure as $\t b$. We refrain from considering this level of generality to keep the notations to a reasonable level. We use the notation $\alpha_t=(\alpha^1_t,\cdots,\alpha^N_t)$ for the players strategies. Notice that the dynamics \reff{fo:firstd} can be rewritten in the form:
\begin{equation}
\label{fo:generald}
dX^i_t=b(t,X^i_t,\overline{\mu}^N_t,\alpha^i_t)dt + \sigma dW^i_t, 
\end{equation}
if the function $b$ of time, a private state, a probability distribution on private states, and a control, is defined by
\begin{equation}
\label{fo:order1}
b(t,x,\mu,\alpha)=\int_{\RR} \tilde b(t,x,x',\alpha)\, d\mu(x')
\end{equation}
and the measure $\overline\mu^N_t$ is defined as the empirical distribution of the private states, i.e.
\begin{equation}
\label{fo:empirical_mu}
\overline{\mu}^N_t=\frac1N\sum_{j=1}^N\delta_{X^j_t}.
\end{equation}
Interactions given by functions of the form \reff{fo:order1} will be called \emph{linear} or \emph{of order $1$}. 
We could imagine that the drift of \reff{fo:firstd} giving the interaction between the private states is of the form
$$
\frac{1}{N^2}\sum_{j,k=1}^N\tilde b(t,X^i_t,X^j_t,X^k_t,\alpha^i_t)
$$
which could be rewritten, still with $\mu=\o\mu^N_t$, \ in the form
\begin{equation}
\label{fo:order2}
b(t,x,\mu,\alpha)=\int_{\RR} \tilde b(t,x,x',x'',\alpha)\, d\mu(x')d\mu(x'').
\end{equation}
Interactions of this form will be called \emph{quadratic} or \emph{of order $2$}. Clearly, one can extend this definition
to interactions of all orders, and more generally, we will say that the interaction is \emph{fully nonlinear}
if it is given by a drift of the form $b(t,X^i_t,\o\mu^N_t,\alpha_t)$ for a general function $b$ defined on
$[0,T]\times\RR\times\cP_{1}(\RR)\times A$ where $\cP_{1}(\RR)$ denotes the set of probability measures on the real line and $A$ the space
in which, at each time $t\in[0,T]$, the controls can be chosen by the individual players.
In general, we say that the game involves interactions of the \emph{mean field} type if the coefficients of the stochastic differential equation giving the dynamics of a private state depend upon the other private states exclusively through the empirical distribution of these private states -- in other words if the interaction is fully nonlinear in the sense just defined -- and if the running and terminal cost functions have the same structure. 

\begin{remark}
As in the case of several of the models discussed later on, the dependence upon the empirical distribution $\o\mu^N_t$ of the private states can degenerate to a dependence upon some moments of this distribution. To be more specific, we can have
\begin{equation}
\label{fo:scalar}
b(t,x,\mu,\alpha)= \tilde b(t,x,\langle\varphi,\mu\rangle,\alpha)
\end{equation}
for some scalar function $\varphi$ of the private states, where we use the duality notation
$
\langle\varphi,\mu\rangle=\int_{\RR}\varphi(x')d\mu(x')
$
for the integral of a function with respect to a measure.
In such a case, we shall say that the interaction is \emph{scalar}.
\end{remark}

\vskip 6pt
To summarize the problem at hand, the game consists in minimizing simultaneously costs of the form
\begin{equation}
\label{fo:Ji}
J^i(\u \alpha)=\EE\left[\int_0^Tf(t,X^i_t,\overline{\mu}^N_t,\alpha^i_t)dt + g(X_T,\overline{\mu}^N_T)\right], \qquad 
i = 1, \cdots, N,
\end{equation}
under constraints of the form
\begin{equation}
\label{fo:dXit}
dX^i_t=b(t,X^i_t,\overline{\mu}^N_t,\alpha^i_t)dt + \sigma dW^i_t, \qquad 0 \leq t \leq T,
\end{equation}
where the $\u W^i$ are independent standard Wiener processes.
(In the whole paper, we use the underlined notation $\u Z$ to denote a stochastic $(Z_{t})_{t \in I}$ indexed by time $t$ in some interval $I$.) 
Note that for symmetry reasons, we choose the running and terminal cost functions $f$ and $g$ to be the same for all the players.
Our goal is to search for equilibriums of such a stochastic differential game.

\subsection{Optimization Problems for Individual Players} 
Given that the problem is intractable in most cases (see nevertheless the discussion of the linear quadratic case in Sections \ref{se:LQMFG}
and \ref{se:LQMKV} and the description of general results toward a systematic treatment in 
Section \ref{se:solvability}), we try to identify realistic models for which  approximate equilibriums and optimal strategies can be identified and computed.

For the sake of definiteness, we restrict ourselves to equilibriums given by Markovian strategies in closed loop feedback form
$$
\alpha_t=(\phi^1(t,X_t),\cdots,\phi^N(t,X_t)), \qquad\qquad 0\le t\le T,
$$
for some deterministic functions $\phi^1$, $\cdots$, $\phi^N$ of time and the state of the system.
Further, we assume that the strategies are of distributed type in the sense that the function $\phi^i$ giving the strategy of player $i$ depends upon the state $X_t$ of the system only through the private state $X_t^i$ of player $i$. In other words, we request that:
$$
\alpha^i_t=\phi^i(t,X^i_t),\qquad\qquad i=1,\cdots,N.
$$
Moreover, given the symmetry of the set-up, we restrict our search for equilibriums to situations in which all the players use the same feedback strategy function, i.e.
$$
\phi^1(t,\,\cdot\,)=\cdots=\phi^N(t,\,\cdot\,)=\phi(t,\,\cdot\,),\qquad\qquad 0\le t\le T,
$$
for some common deterministic function $\phi$.
The search for an equilibrium with a fixed finite number of players involves \emph{optimization}, and the exact form of the optimization depends upon the notion of equilibrium we are interested in. In any case, we hope that, in the large game regime (i.e. when we let the number of players go to $\infty$), some of the features of the problem will be streamlined and the optimization, if feasible, will provide a way to construct
approximate equilibriums for the game with a finite (though large) number of players. We show that, depending upon the nature of the equilibrium we are aiming for, the effective problem appearing in the limit can be of one form or another.

\subsection{Search for Nash Equilibriums: Optimizing First} 
If we search for a Nash equilibrium when the number $N$ of players is fixed, each player $i$ assumes that the other players have already chosen their strategies, say $\alpha^{1*}_t=\phi^{1*}(t,X^1_t)$, $\cdots$ , $\alpha^{i-1*}=\phi^{i-1*}(t,X^{i-1}_t)$, $\alpha^{i+1*}_t=\phi^{i+1*}(t,X^{i+1}_t)$, $\cdots$ , $\alpha^{N*}_t=\phi^{N*}(t,X^N_t)$, and under this assumption, solves the optimization problem:
\begin{equation}
\label{fo:NashN}
\phi^{i*}=\text{arg} \min_\phi \EE\left[\int_0^Tf(t,X^i_t,\overline{\mu}^N_t,\phi(t,X^i_t))dt + g(X_T^i,\overline{\mu}^N_T)\right],
\end{equation}
when the dynamics of his own private state is controlled by the feedback control $\alpha^i_t=\phi(t,X^i_t)$ while the dynamics of the private states of the other players are controlled by the feedbacks $\alpha^j_t=\phi^{j*}(t,X^j_t)$ for $j=1,\cdots,i-1,i+1,\cdots, N$. So when writing an equation for a critical point (typically a \emph{first order condition}), we apply an infinitesimal perturbation to $\phi^*$ without perturbing any of the other $\phi^{j*}$ for $j\ne i$, and look for conditions under which this deviation for player $i$ does not make him better off when the other players do not deviate from the strategies they use. 
Also, recall that we are willing to restrict ourselves to optimal strategies of the form $\phi^{1*}=\phi^{2*}=\cdots=\phi^{N*}$ because of our symmetry assumption. As a result of this \emph{exchangeability} property of the system, even though a perturbation of $\phi^*$ could in principle affect the dynamics of all the private states by changing the empirical distribution of these private states, however, especially when the number of players $N$ is large, a form of the law of large numbers should make this empirical measure $\o\mu^N_t$ quite insensitive
to small perturbations of $\phi^*$. So for all practical purposes,  the optimization problem \reff{fo:NashN} can be solved  as  (or at least its solution can be approximated by) a standard stochastic control problem once the family  $(\overline{\mu}^N_t)_{0\le t\le T}$ of probability measures is fixed. 

\vskip 2pt
So in this approach, the search for an approximate Nash equilibrium is based on the following strategy: first one fixes a family $(\mu_t)_{0\le t\le T}$ of probability measures, next, one solves the standard stochastic control problem (parameterized by the choice of the family $(\mu_t)_{0\le t\le T}$):
\begin{equation}
\label{fo:MFGoptimization}
\phi^{*}=\text{arg} \min_\phi \EE\left[\int_0^Tf(t,X_t,\mu_t,\phi(t,X_t))dt + g(X_T,\mu_T)\right]
\end{equation}
subject to the dynamic constraint
\begin{equation}
\label{fo:MFGconstraint}
dX_t=b(t,X_t,\mu_t,\phi(t,X_t))dt + \sigma d\tilde W_t,
\end{equation}
for some Wiener process $(\tilde W_t)_{0\le t\le T}$. Once this standard stochastic control problem is solved, the remaining issue is the choice of the flow
$(\mu_t)_{0\le t\le T}$ of probability measures. This is where the limit $N\to\infty$ comes to the rescue.  Indeed, the asymptotic independence and the form of law of large numbers provided by the theory of propagation of chaos (see for example \cite{Sznitman}) tell us that, in the limit $N\to\infty$, the empirical measure $\overline{\mu}^N_t$ should coincide with the statistical distribution of $X_t$.  So, once an optimum $\phi^*$ is found for each choice of the family $(\mu_t)_{0\le t\le T}$, then the family of measures is determined (typically by a fixed point argument) so that, at each time $t$, the statistical distribution of the solution $X_t$ of \reff{fo:MFGconstraint} is exactly $\mu_t$. 
\vskip 2pt
The role of the limit $N\to\infty$ (large number of players) is to guarantee the stability of the empirical measure $\o\mu^N_t$ when a single player perturbs his strategy while the other players keep theirs unchanged, and the fact that this stable distribution has to be the common statistical distribution of all the private states $X^i_t$.
Performing the optimization over $\phi$ when the family $(\mu_t)_{0\le t\le T}$ of probability measures is kept fixed is the proper way to implement the notion of Nash equilibrium whereby each player is not better off if he deviates from his strategy while all the other players keep theirs untouched, as implied by the lack of change in  
$(\mu_t)_{0\le t\le T}$ and hence $(\overline\mu^N_t)_{0\le t\le T}$. 
\vskip 2pt

From a mathematical point of view, the limit $N \to \infty$ can be justified rigorously in the following sense: under suitable conditions, if there exists a family of measures $(\mu_{t})_{0 \leq t \leq T}$ together with an optimally controlled process in 
(\ref{fo:MFGoptimization}--\ref{fo:MFGconstraint})
with $(\mu_{t})_{0 \leq t \leq T}$ as marginal distributions exactly, then the corresponding optimal feedback $\phi^*$ is proven to provide an $\varepsilon$-Nash equilibrium to the $N$-player-game \eqref{fo:NashN} for $N$ large. (See \cite{Cardaliaguet,CarmonaDelarue2}.) 
The argument relies on standard theory of propagation of chaos for McKean--Vlasov diffusion processes. (See \cite{Sznitman}.) 
As a by-product, the empirical measure of the $N$-player system, when driven by $\phi^*$, does indeed converge toward 
$(\mu_t)_{0\le t\le T}$. In such cases, this makes rigorous the MFG approach for searching for an approximate Nash equilibrium to the $N$-player-game. Still the converse limit property seems to be much more involved and remains widely open: do the optimal states of the $N$-player-game (if exist) converge to some optimal states of 
(\ref{fo:MFGoptimization}--\ref{fo:MFGconstraint}) as $N$ tends to $\infty$?  
\vskip 4pt

Assuming that the family $\underline\mu=(\mu_t)_{0\le t\le T}$ of probability distributions has been frozen,
the Hamiltonian of the stochastic control problem is
\begin{equation}
\label{fo:hmu}
H^{\mu_{t}}(t,x,y,\alpha)=y b(t,x,\mu_t,\alpha) + f(t,x,\mu_t,\alpha).
\end{equation}
Recall that we are limiting ourselves to the case of constant and identical volatilities $\sigma$ for the sake of simplicity.
 In all the examples considered in this paper, there exists a regular function $\hat \alpha$ satisfying:
\begin{equation}
\label{fo:alphahat}
\hat \alpha=\hat \alpha^{\mu_{t}}:[0,T]\times \RR\times\RR \ni(t,x,y) \hookrightarrow \hat \alpha(t,x,y)\in\text{arg}\min_{\alpha\in A}H^{\mu_{t}}(t,x,y,\alpha).
\end{equation}
We denote by $\cH^{\mu_{t}}(t,x,y)$ this infimum:
\begin{equation}\label{hamiltonian}
\cH^{\mu_{t}}(t,x,y)=\inf_{\alpha\in A} H^{\mu_{t}}(t,x,y,\alpha).
\end{equation}
In this paper, we solve stochastic control problems using the probabilistic approach based on the Pontryagin minimum principle and solving the adjoint forward-backward stochastic differential equations. For the sake of completeness, we review the Lasry--Lions' approach using the Hamilton--Jacobi--Bellman (HJB for short) equation leading to the solution of a forward-backward system of nonlinear Partial Differential Equations (PDEs for short).

\subsubsection*{HJB / PDE approach}
Since ${\u\mu}$ is frozen, the stochastic control problem is Markovian and we can introduce the HJB value function:
\begin{eqnarray}
v(t,x) = \inf _{\alpha\in \cA_t} \EE \left[ \int_t^T    f(s,X_s,\mu_s,\alpha_s)ds+g(X_T,\mu_T)| X_t=x\right]. \label{eq:16:12:1}
\end{eqnarray}
where $\cA_t$ denotes the set of admissible controls over the interval $[t,T]$. We expect that the HJB value function $v$ is the solution in some sense (most likely in the viscosity sense only) of the Hamilton--Jacobi--Bellman (HJB) equation (see \cite{FlemingSoner})
\begin{equation}\label{hjb}
\partial_tv+\frac{\sigma^2}{2}\partial^2_{xx} v +\cH^{\mu_{t}}(t,x, \partial_x  v(t,x)) = 0, \qquad 
(t,x) \in [0,T] \times \RR,
\end{equation}
with terminal condition $v(T,x)=g(x,\mu_T)$, $x \in \RR$.

\vskip 2pt
In the end, we would like $\underline\mu=(\mu_t)_{0\le t\le T}$ to be the flow of marginal distributions of the optimally controlled private state. Clearly, this sounds like a circular argument at this stage since this stochastic differential equation actually depends upon ${\u\mu}$, and it seems that only a fixed point argument can resolve such a quandary. In any case, the flow of statistical distributions should satisfy Kolmogorov's equation. In other words, if we use the notation
$$
\beta(t,x)= b(t,x,\mu_t,\phi(t,x))
$$
where $\phi$ is the optimal feedback control (think of $\phi(t,x)=\hat \alpha^{\mu_{t}}(t,x,\partial_xv(t,x))$), then the flow $(\nu_t)_{0\le t\le T}$ of measures should be given by $\nu_t=\cL(X_t)$ and satisfy Kolmogorov's equation
\begin{equation}
\label{fo:kolmo}
\partial_t\nu-\frac{\sigma^2}{2}\partial^2_{xx} \nu -\text{div}(\beta(t,x)\nu) = 0, \qquad 
(t,x) \in [0,T] \times \RR,
\end{equation}
with initial condition $\nu_0=\mu_0$. This PDE can be given a rigorous meaning in the sense of distributions.  When $\nu_t$ has a smooth density, integration by part can be used to turn the formal PDE \reff{fo:kolmo} for $\nu$ into a classical PDE for the density of $\nu$.

\vskip 2pt
Setting  $\nu_t=\cL(X_t)=\mu_t$ when $(X_t)_{0\le t\le T}$ is the diffusion optimally controlled by the feedback function $\phi$
gives a system of coupled nonlinear forward-backward PDEs \pref{hjb}--\pref{fo:kolmo}, called the MFG PDE system. See \cite{MFG1,MFG2,MFG3,MFG4,Lions,GueantLasryLions.pplnm}.

\subsubsection*{Stochastic Maximum Principle / FBSDE Approach}
Given the frozen flow of measures $\u\mu$, the stochastic optimization problem is a standard stochastic control problem and as such, its solution can be approached via the stochastic Pontryagin principle. For each open-loop adapted control $\u\alpha=(\alpha_t)_{0\le t \le T}$, we denote by $\u X^{\u\alpha}=(X^{\u\alpha}_t)_{0\le t\le T}$ the associated state; any solution $(Y_t,Z_t)_{0\le t\le T}$ of the BSDE
\begin{equation}
\label{fo:mfgadjointeqn}
dY_t=- \partial_{x} H^{\mu_{t}}\bigl(t,X_{t},Y_{t},\alpha_t\bigr) dt + Z_{t} dW_{t},
\qquad t \in [0,T] \ ; \qquad Y_{T} = \partial_{x} g(X_{T},\mu_{T}),
\end{equation}
is called a set of adjoint processes and the BSDE (\ref{fo:mfgadjointeqn}) is called the adjoint equation. 
The necessary condition of the stochastic Pontryagin principle says that, whenever $\u X^{\u \alpha}$ is an optimal state of the optimization problem, it must hold $H^{\mu_{t}}(t,X_{t},Y_{t},\alpha_{t}) = {\mathcal H}^{\mu_{t}}(t,X_{t},Y_{t})$, 
$t \in [0,T]$. (See \cite{YongZhou}.)
Conversely, when the Hamiltonian $H^{\mu_{t}}$ is convex with respect to  the variables 
$(x,\alpha)$ and the terminal cost $g$ is convex with respect to  the variable $x$, the forward component of any solution to the FBSDE
\begin{equation*}
\begin{cases}
&dX_{t} = b\bigl(t,X_{t},\mu_{t},\hat{\alpha}^{\mu_{t}}(t,X_{t},Y_{t})\bigr) dt + \sigma dW_{t}, 
\\
&dY_{t} = - \partial_{x} H^{\mu_{t}}\bigl(t,X_{t},Y_{t},\hat{\alpha}^{\mu_{t}}(t,X_{t},Y_{t})\bigr) dt + Z_{t} dW_{t},
\end{cases}
\end{equation*}
with the right initial condition for $X_{0}$ and the terminal condition $Y_{T} =  \partial_{x}g(X_{T},\mu_{T})$, is an optimally controlled path for \eqref{fo:MFGoptimization}. 
In particular, if we want to include the matching condition of the MFG approach as outlined earlier, namely if we want to enforce $\mu_t=\cL(X_t)$, the 
above stochastic forward-backward system turns into
\begin{equation}
\label{eq:16:12:2}
\begin{split}
&dX_{t} = b\bigl(t,X_{t},\cL(X_{t}),\hat{\alpha}^{\cL(X_{t})}(t,X_{t},Y_{t})\bigr) dt + \sigma dW_{t}, 
\\
&dY_{t} = - \partial_{x} H^{\cL(X_{t})}\bigl(t,X_{t},Y_{t},\hat{\alpha}^{\cL(X_{t})}(t,X_{t},Y_{t})\bigr) dt + Z_{t} dW_{t},
\end{split}
\end{equation}
with the right initial condition for $X_{0}$ and with the terminal condition $Y_{T} =  \partial_{x}g(X_{T},\cL(X_{T}))$. This FBSDE is of the McKean--Vlasov type as formally introduced in Subsection \ref{subsec:generalMKV}. If convexity (as described above) holds, 
it characterizes the optimal states of the MFG problem. If convexity fails, it provides a necessary condition only. 

\subsection{Search for Cooperative Equilibriums: Taking the Limit $N\to\infty$ First} 
\label{subse:COOP}
As before, we assume that all the players use the same feedback function $\phi$, but the notion of optimality is now defined in terms of what it takes for a feedback function $\phi^*$ to be critical. In the search for a first order condition, for each player $i$, when $\phi^*$ is perturbed into $\phi$, we now assume that all the players  $j\ne i$ also use the perturbed feedback function. In other words, they use the controls $\alpha^j_t=\phi(t,X^j_t)$ for $j=1,\cdots,i-1,i+1,\cdots, N$. In this case, the perturbation of $\phi^*$ has a significant impact on the empirical distribution $\o\mu^N_t$, and the optimization cannot be performed after the latter is frozen, as in the case of the search for a Nash equilibrium. In the present situation, the simplification afforded by taking the limit $N\to\infty$ is required before we perform the optimization.

\vskip 4pt
Assuming that the feedback function $\phi$ is fixed, at least temporarily, the theory of propagation of chaos states that, if we consider the solution $X^N_t=(X^{N,1}_t,\cdots, X^{N,N}_t)$ of the system of $N$ stochastic differential equations \reff{fo:dXit} with $\alpha^i_t=\phi(t,X^{N,i}_t)$, then in the limit $N\to \infty$, for any fixed integer $k$,  the joint distribution of the $k$-dimensional process $\{(X^{N,1}_t,\cdots,X^{N,k}_t)\}_{0\le t\le T}$ converges to a product distribution (in other words the $k$ processes $(X^{N,i}_t)_{0\le t\le T}$ for $i=1,\cdots,k$ become independent in the limit) and the distribution of each single marginal process converges toward the distribution of the unique solution $\u X=(X_t)_{0\le t\le T}$ of the McKean--Vlasov evolution equation
\begin{equation}
\label{fo:mkv}
d X_t=b(t,X_t,\cL(X_t),\phi(t,X_t))dt + \sigma d\tilde W_t
\end{equation}
where $(\tilde W_t)_{0\le t\le T}$ is a standard Wiener process. So if the common feedback control function $\phi$ is fixed, in the limit $N\to\infty$, the private states of the players become independent of each other, and for each given $i$, the distribution of the private state process $(X^{N,i}_t)_{0\le t\le T}$ evolving according to \reff{fo:dXit} converges toward the distribution of the solution of \reff{fo:mkv}. So if we optimize after taking the limit $N\to\infty$, i.e. assuming that the limit has already been taken, the objective of each player becomes the minimization of the functional
\begin{equation}
\label{fo:mkvobjective}
J(\phi)=\EE\left[\int_0^Tf(t,X_t,\cL(X_t),\phi(t,X_t))dt + g(X_T,\cL(X_T))\right]
\end{equation}
over a class of admissible feedback controls $\phi$. Minimization of \eqref{fo:mkvobjective} over $\phi$ under the dynamical constraint \eqref{fo:mkv} is a form of optimal stochastic control where the controls are in closed loop feedback form. More generally, such a problem can be stated as in \eqref{fo:MFGoptimization} for open loop controls $\u\alpha=(\alpha_t)_{0\le t \le T}$ adapted to any specific information structure: 
\begin{eqnarray}
\label{fo:mkvcontrolpb}
\u\alpha^*&=&\text{arg}\min_{\alpha}\EE\left[\int_0^Tf(t,X_t,\cL(X_t),\alpha_t)dt + g(X_T,\cL(X_T))\right]\nonumber\\
\text{subject to}&&\\  
&&\phantom{????}dX_t=b(t,X_t,\cL(X_t),\alpha_t)dt + \sigma dW_t, \qquad t \in [0,T].\nonumber
\end{eqnarray}
Naturally, we call this problem the \emph{optimal control of the stochastic McKean-Vlasov dynamics}. 
Here as well, the limit $N \rightarrow \infty$ can be justified rigorously: under suitable conditions, optimal feedback controls in 
\eqref{fo:mkvcontrolpb} are proven to be $\varepsilon$-optimal controls for the $N$-player-system when $N$ is large, provided that the rule we just prescribed above is indeed in force, that is: any perturbation of the feedback is felt by all the players in a similar way. (See \cite{CarmonaDelarue4}.)

\vskip 4pt
Standard techniques from stochastic control of Markovian systems cannot be used for this type of stochastic differential equations
and a solution to this problem is not known in general, even though solutions based on appropriate analogs of the Pontryagin minimum principle have been sought for by several authors. 
We first review the results of \cite{AndersonDjehiche} which are the only published ones which we know of. They only concern scalar interactions  for which the dependence upon the probability measure of the drift and cost functions are of the form
\begin{equation*}
b(t,x,\mu,\alpha)= b(t,x,\langle\psi,\mu\rangle,\alpha), \quad 
f(t,x,\mu,\alpha) = f(t,x,\langle\gamma,\mu\rangle,\alpha), \quad
g(x,\mu)=g(x,\langle\zeta,\mu\rangle)
\end{equation*}
for some functions $\psi$, $\gamma$, and $\zeta$. Recall that we use the duality notation $\langle\varphi,\mu\rangle$ to denote the integral of the function $\varphi$ with respect to the measure $\mu$. Notice that we use the same notations $b$, $f$ and $g$ for functions where the variable $\mu$, which was a measure, is replaced by a numeric variable.
As we explained earlier, this setting could be sufficient for the analysis of most of the linear quadratic models we consider in Sections \ref{se:LQMFG} and \ref{se:LQMKV}. We review the more general results \cite{CarmonaDelarue2} of the first two named authors in Section \ref{se:solvability}. 

The major difference with the classical case is in the form of the adjoint equation. Given a control process $\u\alpha=(\alpha_t)_{0\le t \le T}$ and a process $\u X=(X_t)_{0\le t\le T}$ satisfying
\begin{equation}
\label{fo:generalmkv}
dX_t=b(t,X_t,\cL(X_t),\alpha_t)dt + \sigma dW_t, \qquad 0 \leq t \leq T,
\end{equation}
a pair of processes $\u Y=(Y_t)_{0\le t\le T}$  and $\u Z=(Z_t)_{0\le t\le T}$  is said to form a pair of adjoint processes if they satisfy:
\begin{equation}
\label{fo:adjointmkv:0}
\begin{cases}
dY_t&=-[\partial_xb(t,X_t,\EE\{\psi(X_t)\},\alpha_t)Y_t+ \partial_xf(t,X_t,\EE\{\gamma(X_t)\},\alpha_t)]dt+Z_td  W_t\\
&\hspace{-25pt}-[\EE\{\partial_{x'}b(t,X_t,\EE\{\psi(X_t)\},\alpha_t)Y_t\}\partial_x\psi(X_t) +\EE\{ \partial_{x'}f(t,X_t,\EE\{\gamma(X_t)\},\alpha_t)\}\partial_x\gamma(X_t)]dt
\\
Y_T&=\partial_x g(X_T,\EE\{\zeta(X_T)\}) + \EE\{ \partial_{x'}g(X_T,\EE\{\zeta(X_T)\})\}\partial_x\zeta(X_T).
\end{cases}
\end{equation}
This BSDE, which is of the McKean--Vlasov type as its coefficients depend upon the distribution of the solution, is called the adjoint equation. It provides a necessary condition for the optimal states of the MKV optimal control problem, that is, for any optimally controlled 
path, it must hold:
\begin{equation}
\label{fo:isaacsmkv}
\tilde H(t,X_t,Y_t,\alpha_t) =\inf_{\alpha\in A}\tilde H(t,X_t,Y_t,\alpha)
\end{equation}
for all $t\in[0,T]$. The Hamiltonian $\tilde H$ appearing in \reff{fo:isaacsmkv} is defined, for any random variable $\xi$, as:
\begin{equation}
\label{fo:hamiltonianmkv:2}
\tilde H(t,\xi,y,\alpha)=yb(t,\xi,\EE\{\psi(\xi)\},\alpha)  + f(t,\xi,\EE\{\gamma(\xi)\},\alpha).
\end{equation}
This result is proven in  \cite{AndersonDjehiche} (with the assumption that $A$ is convex and that $\tilde{H}$ is convex in the variable $\alpha$). The sufficiency condition is given by (see \cite{AndersonDjehiche} as well):
\begin{theorem}
\label{th:pontryaginmkv}
Assume that
\begin{enumerate}\itemsep=-2pt
\item[(A1)] $g$ is convex in $(x,x')$;
\item[(A2)] the partial derivatives $\partial_{x'}f$ and $\partial_{x'}g$ are non-negative;
\item[(A3)] the Hamiltonian $H$ is convex in $(x,x_{1}',x_{2}', \alpha)$;
\item[(A4)] the function $\psi$ is affine and the functions $\gamma$ and $\zeta$ are convex;
\end{enumerate}
where the Hamiltonian function $H$ is defined as:
\begin{equation}
\label{fo:hamiltonianmkv}
H(t,x,x_{1}',x_{2}',y,\alpha)=yb(t,x,x_{1}',\alpha)  + f(t,x,x_{2}',\alpha).
\end{equation}
If $\u X=(X_{t})_{0 \leq t \leq T}$ satisfies \eqref{fo:generalmkv}
for some control process $\u\alpha=(\alpha_t)_{0\le t \le T}$, if $\u Y=(Y_t)_{0\le t\le T}$  and $\u Z=(Z_t)_{0\le t\le T}$ form a pair of  adjoint processes
and if \eqref{fo:isaacsmkv} holds for all $t\in[0,T]$ (almost surely), then the control $\underline{\alpha}
= (\alpha_{t})_{0 \leq t \leq T}$ is optimal. 
\end{theorem}

\subsection{Summary}
The dichotomy in purpose suggested by the two notions of equilibrium discussed above points to two different paths to go from the North East corner
to the South West corner of the following diagram, and the thrust of the paper is to provide insight, and demonstrate by examples that which path one chooses has consequences on the properties of the equilibrium, in other words, the diagram is not \emph{commutative}.
$$
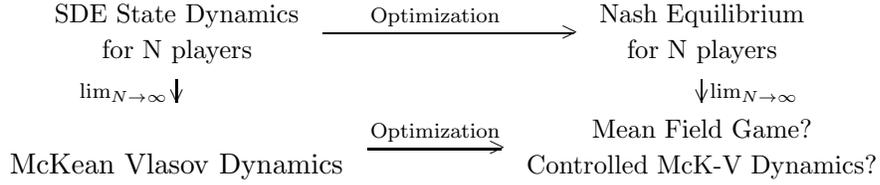
\begin{diagram} 
{\begin{array}{c}\text{\small SDE State Dynamics}\\ \text{\small for N players}\\  \end{array}}& \rTo^{\text{Optimization}}& {\begin{array}{c}\text{\small Nash Equilibrium}\\ \text{\small for N players}\\ 
\end{array}}\\
\dTo^{\lim_{N\to\infty}} &	&\dTo_{\lim_{N\to\infty}}\\
{\begin{array}{c}\text{\small }\\ \text{McKean Vlasov Dynamics}\end{array}}	&\rTo^{\text{Optimization}}	&\begin{array}{c}\text{\small Mean Field Game?} \\ \text{\small Controlled McK-V Dynamics?} \end{array} 
\end{diagram}
$$
It is important to emphasize one more time what we mean by the limit $N\to\infty$. We want to identify properties of the limit which, when re-injected into the game with finitely many players, give an approximate solution to the problem we are unable to solve directly for the stochastic game with $N$ players.

\subsubsection*{First Discussion of the Differences}
In the present context where we assume that the volatility $\sigma$ is constant, the general form of MKV dynamics is given by the solution of non-standard stochastic differential equations in which the distribution of the solution appears in the coefficients:
\begin{equation}
\label{X_MKV} 
dX_t= b(t,X_t,\cL(X_t),\alpha_t) dt+\sigma dW_t, \; \qquad X_0=x, 
\end{equation}
where $\sigma$ is a positive constant. For a given admissible control $\u\alpha=(\alpha_t)_{0\le t\le T}$ we write $\u X^{\u \alpha}$  for the unique solution to \reff{X_MKV}, which exists under the usual growth and Lipchitz conditions on the function $b$, see \cite{Sznitman} for example. The problem is to optimally control this process so as to minimize the expectation:
\begin{equation}
\label{J_MKV}
J_{\rm MKV}(\u\alpha)= \E \left[ \int_0^T    f(t,X_t^{\u\alpha},\cL(X^{\u\alpha}_t),\alpha_t)dt+g(X_T^{\u\alpha},\cL(X_T^{\u\alpha}))\right]. 
\end{equation}
On the other hand, the dynamics arising in the MFG approach are given by the solution of a standard stochastic differential equation:
\begin{equation}
\label{X_MFG} dX_t= b(t,X_t,\mu_t,\alpha_t) dt+\sigma dW_t, \qquad X_0=x, 
\end{equation}
where  $\u\mu=(\mu_t)_{0\le t\le T}$ is a deterministic function with values in the space of probability measures, which can be understood as a (temporary) proxy or candidate for the anticipated statistical distribution of the random variables $X_t$. The expectation that has to be minimized is now:
\begin{equation}
\label{J_MFG}
J_{\rm MFG}(\u \alpha)= \E \left[ \int_0^T    f(t,X^{\u\alpha}_t,\mu_t,\al_t)dt+g(X^{\u\alpha}_T,\mu_T)\right],
\end{equation}
and an equilibrium takes place if, for an optimal control ${\u\alpha}$, the anticipated distribution $\mu_t$ actually coincides with $\cL (X^{\u\alpha}_t)$ for every $t\in[0,T]$. While very similar, the two problems differ by a very important point: the timing of  the optimization. 

In the MKV control problem (\ref{X_MKV})-(\ref{J_MKV}), at each time $t$, the statistical distribution $\mu_t$ is matched to the distribution of the state, and once we have $\mu_t=\cL(X_t)$, then we do optimize the objective function. On the other hand,  for the MFG problem (\ref{X_MFG})--(\ref{J_MFG}), the optimization is performed for each fixed family $\u\mu=(\mu_t)_{0\le t\le T}$ of probability distributions, and once the optimization is performed, one matches the distributions $\mu_t$ to the laws of the optimally controlled process.

\section{Mean Field Linear Quadratic (LQ) Games}
\label{se:LQMFG}
This section is devoted to a class of models for which we can push the analysis further, to the point of deriving explicit formulas in some cases. We choose coefficients of the linear-quadratic type, that is
\begin{equation*}
\begin{split}
&b(t,x,\mu,\alpha) = a_tx+\o a_t\o\mu+ b_t\alpha+\beta_t,
\\
&f(t,x,\mu,\alpha) = \frac12 n_{t}\alpha^2 + \frac12
( m_t x+\o m_t \overline{\mu}_t)^2, \quad g(x,\mu) = \frac12 (q x+\o q \, \o\mu)^2
\end{split}
\end{equation*}
where $a_{t}$, $\o a_{t}$, $b_{t}$, $\beta_{t}$, $m_{t}$, $\o m_{t}$ and $n_{t}$
are deterministic continuous functions of $t\in [0,T]$, and $q$ and $\o q$ are deterministic. Recall that we use the notation $\o\mu$ for the mean of the measure $\mu$, i.e. $\o\mu=\int xd\mu(x)$.

\subsection{Solving the $N$-player Game}
\label{subse:solvNPG}
Under specific assumptions on the coefficients, existence of open loop Nash equilibriums for Linear Quadratic (LQ for short) stochastic differential games of the type considered in this section can be proven using the stochastic Pontryagin 
principle approach. See for example \cite{Hamadene}. We re-derive these results in the case of mean field models  considered in this paper. 

\vskip 1pt
As in our introductory discussion, we assume that individual private states are one dimensional.
This assumption is not essential, its goal is only to keep the notation to a reasonable level of complexity. See nevertheless Remark \ref{re:ricatti} below. Moreover, we also assume that the actions of the individual players are also one-dimensional (i.e. $A=\RR$). Using lower cases for all the examples tackled below, the dynamics of the private states rewrite:
$$
dx^i_t=[a_tx^i_t +\o a_t\overline x_t +b_t\alpha^i_t+\beta_t]dt + \sigma dW^i_t,\qquad i=1,\cdots,N, \qquad t \in [0,T]; \qquad 
x_{0}^i = x^{(0)},
$$
where $x^{(0)}$ is a given initial condition in $\RR$ and $\overline{x}_{t}$ stands for the 
empirical mean:
\begin{equation*}
\overline x_{t} = \frac{1}{N} \sum_{i=1}^N x_{t}^i.
\end{equation*}
Similarly, the individual costs become:
$$
J^i(\u \alpha)=\EE\left\{\int_0^T\frac12 \bigl[n_t (\alpha_t^i)^2+(m_t x_t^i+{\o m}_t\o x_t)^2\bigr]dt + \frac12(qx_T^i +{\o q} \,\o x_T)^2\right\}, \qquad i = 1, \cdots, N.
$$
We will assume throughout the section that 
$\inf_{t \in [0,T]} b_t>0$ and $\inf_{t \in [0,T]} n_t>0$.

\subsection{Implementing the MFG Approach}
As explained earlier, the first step of the MFG approach is to fix a flow $\u\mu=(\mu_t)_{0\le t\le T}$ of probability measures in lieu of the empirical measures of the players' private states, and solve the resulting control problem for one single player. Since the empirical measure  of the players' private states enters only the state equations and the cost functions through its mean, it is easier to choose a real valued deterministic function $(\o \mu_t)_{0\le t\le T}$ (which we should denote  $\o{\u \mu}$ in order to be consistent with the convention used so far in our notation system)
as a proxy for the mean of the empirical distribution of the private states. Then the individual stochastic control problem is to minimize
$$
J(\u \alpha)=\EE\left\{\int_0^T[\frac{n_t}2\alpha_t^2+ \frac12 (m_tx_t+{\o m}_t\o \mu_t)^2]dt + \frac12(qx_T+{\o q} \, \o \mu_T)^2\right\}
$$
subject to the dynamical constraint
$$
dx_t=[a_tx_t +{\o a}_t\o \mu_t +b_t\alpha_t+\beta_t]dt + \sigma dW_t, \qquad t \in [0,T],
$$
with $x_{0}=x^{(0)}$ as initial condition.
Given that the deterministic function $t\hookrightarrow \o \mu_t$ is assumed fixed, the Hamiltonian of the stochastic control problem is equal to
\begin{equation}
\label{fo:hamiltonianmfg}
H^{\o{\mu}_{t}}(t,x,y,\alpha)=y[a_t x +{\o a}_t
 \o\mu_t + b_t \alpha +\beta_t]+\frac12(m_t x +{\o m}_t\o \mu_t)^2 +\frac12 n_t \alpha^2. 
\end{equation}
The optimal control $\hat \alpha$ minimizing the Hamiltonian is given by:
\begin{equation}
\label{fo:alphahat1}
\hat\alpha(t,y)=-\frac{b_t}{n_t}y, \qquad t \in [0,T], \ y \in \RR,
\end{equation}
and the minimum Hamiltonian by
$$
\cH^{\o{\mu}_{t}}(t,x,y)=\inf_{\alpha\in A}H^{\o{\mu}_{t}}(t,x,y,\alpha)=(a_t x+\o a_{t}
\o \mu_t+\beta_t)y +\frac12(m_tx+\o m_{t}  \o \mu_t)^2-\frac12\frac{b^2_t}{n_t}y^2.
$$

\subsubsection*{Stochastic Pontryagin Principle Approach}
In the present situation, the adjoint equation reads:
\begin{equation}
\label{fo:adjointmfg}
dy_t=-[a_ty_t+m_t(m_t x_t +{\o m}_t\o \mu_t)]dt+z_tdW_t,\qquad
y_T= q(q x_T +{\o q} \, \o \mu_T).
\end{equation}
Notice that because of the special form of the coefficients, this BSDE does not involve the control $\hat\alpha_t$ explicitly.
By \eqref{fo:alphahat1}, the optimal control has the form $\hat\alpha_t=-n_{t}^{-1}b_{t} y_{t}$, so that the forward dynamics of the state become:
\begin{equation}
\label{fo:fmfg}
dx_t=[a_tx_t+{\o a}_t\o \mu_t -\frac{b_t^2}{n_t}y_t+\beta_t]dt+\sigma dW_t, \qquad x_{0}=x^{(0)},
\end{equation}
and the problem reduces to the proof of the existence of a solution to the standard FBSDE 
\reff{fo:adjointmfg}--\reff{fo:fmfg}  and the analysis of the properties of such a solution. In particular, we will need to check that this solution is amenable to the construction of a fixed point for $\o{\u\mu}$. Notice that, once this FBSDE is solved, \reff{fo:alphahat1} will provide us with an open loop optimal control for the problem. 

In a first time, we remind the reader of elementary properties of linear FBSDEs. 
We thus streamline the notation and rewrite the FBSDE as
\begin{equation}
\label{fo:xyfbsde}
\begin{cases}
dx_t&=[\fa_tx_t+\fb_ty_t+\fc_t]dt+\sigma dW_t,\qquad x_0=x^{(0)},\\
dy_t&=[\fm_tx_t- \fa_ty_t+ \fd_t]dt+z_tdW_t,\qquad
y_T= \fq x_T +\fr,
\end{cases}
\end{equation}
where in order to better emphasize the linear nature of the FBSDE we set:
\begin{equation*}
\fa_t=a_t, \; \fb_t=- b_t^2/ n_t,\; \fc_t=\beta_t+{\o a}_t\o \mu_t, \; 
\fm_t= - m_t^2, \; \fd_t=- m_t{\o m}_t\o \mu_t, \; \fq=q^2,\; \fr=q{\o q} \, \o \mu_T.
\end{equation*}
Linear FBSDEs of the form \eqref{fo:xyfbsde} have been studied in \cite{Yong1,Yong2}, but for the sake of completeness we construct a solution from scratch recalling some basic facts of the theory. Because of the linearity of the system \reff{fo:xyfbsde}, we expect the FBSDE value function to be affine in the space variable, so we search for deterministic functions $\eta_t$ and $\chi_t$ such that 
\begin{equation}
\label{fo:ansatz:0}
y_t=\eta_tx_t+\chi_t,\qquad  t \in [0,T].
\end{equation}
Computing $dy_t$ from this ansatz using the expression of $dx_t$ from \reff{fo:xyfbsde} we get
$$
dy_t=[(\dot\eta_t+\fa_t\eta_t+\fb_t\eta_t^2)x_t+\dot\chi_t+\fb_t\eta_t\chi_t+\fc_t\eta_t]dt + \sigma\eta_t dW_t, 
\qquad  t \in [0,T],
$$
and identifying term by term with the expression of $dy_t$ given in \reff{fo:xyfbsde} we get:
\begin{equation}
\label{fo:mfgfinal}
\begin{cases}
&\dot\eta_t=-\fb_t\eta_t^2-2\fa_t\eta_t+\fm_t,\qquad \eta_T=\fq\\
&\dot\chi_t+(\fa_t+ \fb_t\eta_t)\chi_t=\fd_t-\fc_t\eta_t,\qquad
\chi_T=\fr\\
&z_t=\sigma \eta_t.
\end{cases}
\end{equation}
The first equation is a Riccati equation. According to the classical theory of Ordinary Differential Equations (ODEs for short), its solution may be obtained by solving the second order linear equation
$$
-\fb_t\ddot \theta_t + [\dot{\fb_t}-2 \fa_t\fb_t]\dot \theta_t + \fm_t\fb_t^2\theta_t=0,
$$
with terminal conditions $\theta_T=1$ and $\dot \theta_T= \fb_T \fq$, and setting
$\eta_t=  (\fb_t \theta_{t})^{-1} \dot \theta_t$,
provided there is a solution $t\hookrightarrow\theta_t$ which does not vanish. 
%
In the framework of control theory, the first equation in 
\eqref{fo:mfgfinal} can be also reformulated as the Riccati equation deriving from a deterministic control problem of linear-quadratic type. Namely, the adjoint system associated with the minimization of the cost function (with the same 
notations as in \eqref{fo:xyfbsde}):  
\begin{equation}
\label{eq:29:12:13:b}
\o J(\u {\o \alpha}) = \frac{1}{2} \fq \xi_{T}^2
+
\int_{0}^T \frac{1}{2} \bigl[ \o \alpha_{t}^2 - \fm_{t} \xi_{t}^2 \bigr] dt
\end{equation}
over the control $(\bar{\alpha}_{t})_{0 \leq t \leq T}$,
subject to  
\begin{equation*}
d \xi_{t}= [ \fa_{t} \xi_{t} - (-\fb_{t})^{1/2} \o \alpha_{t}] dt, \qquad t \in [0,T]; \qquad \xi_{0}=x^{(0)},
\end{equation*}
admits the first equation in \eqref{fo:mfgfinal} as Riccati factorization. Here, $\fq$ and $-\u \fm$ are non-negative so that the 
minimization problem \eqref{eq:29:12:13:b} has a unique optimal path and the Riccati equation is solvable. (See Theorem 37 in \cite{Sontag}.) Therefore, the 
system \eqref{fo:mfgfinal} is solvable. 

As a by-product, we deduce that the forward-backward system \eqref{fo:xyfbsde} is solvable. By the stochastic minimum principle, the solution is even unique. Indeed, the forward-backward system \eqref{fo:xyfbsde} describes the optimal states of the  
stochastic control problem driven by the cost functional 
\begin{equation}
\label{eq:21:6:1}
{\mathfrak J}(\u {\bfm \alpha}) = {\mathbb E} \biggl[ \frac{1}{2} \fq {\bfm \xi}_{T}^2
+ \int_{0}^T \frac{1}{2} \bigl[  {\bfm \alpha}_{t}^2 - \fm_{t} {\bfm \xi}_{t}^2 - 2\fd_{t} {\bfm \xi}_{t}\bigr] dt \biggr]
\end{equation}
depending upon the control $({\bfm \alpha}_{t})_{0 \leq t \leq T}$,
subject to  
\begin{equation}
\label{eq:21:6:2}
d {\bfm \xi}_{t}= [ \fa_{t} {\bfm \xi}_{t} - (-\fb_{t})^{1/2} {\bfm \alpha}_{t} + \fc_{t}] dt + \sigma dW_{t}, \quad t \in [0,T]; 
\qquad {\bfm \xi}_{0}=x^{(0)}. 
\end{equation}
Since the cost coefficients in ${\mathfrak J}$ are convex in $x$ and strictly convex in $\alpha$ and since the drift in 
\eqref{eq:21:6:2} is linear in $x$ and $\alpha$, the stochastic minimum principle says that the optimal control must be unique.  Uniqueness of the solution to \eqref{fo:xyfbsde} follows directly. (See \cite{PengWu}.) Notice  that the convexity of the cost functions holds in the original problem (\ref{fo:adjointmfg}--\ref{fo:fmfg}).
\vspace{4pt}

Going back to the notations in \reff{fo:adjointmfg}--\eqref{fo:fmfg}, we deduce that the system \reff{fo:adjointmfg}--\eqref{fo:fmfg} is always uniquely solvable when $\o \mu$ is given. The point is thus to construct a fixed point for $\o {\u \mu} = (\o \mu_{t})_{0 \leq t \leq T}$.  
If a fixed point does exist, that is if we can solve the system
\reff{fo:adjointmfg}--\eqref{fo:fmfg} with the constraint $\o \mu_{t} = \EE(x_{t})$ for $t \in [0,T]$, then the pair
$(\o \mu_{t}=\EE(x_{t}),\o y_{t} = \EE(y_{t}))_{0 \leq t \leq T}$ solves the deterministic forward-backward system:
\begin{equation}
\label{eq:29:12:10}
\begin{split}
&d\o \mu_t=[(a_t +{\o a}_t) \o \mu_t -\frac{b_t^2}{n_t} \o y_t+ \beta_t]dt, \qquad \o \mu_{0}= x_{0}^{(0)},
\\
&d \o y_{t} = -[a_t \o y_t+m_t(m_t +{\o m}_t) \o \mu_t)] dt,\qquad
\o y_T= q(q +{\o q}) \o \mu_T.
\end{split}
\end{equation}
Conversely, if the system \eqref{eq:29:12:10} is solvable, then $\u {\o \mu}$ satisfies the fixed point condition for the mean-field game of linear-quadratic type since the linear system
\begin{equation}
\label{eq:31:12:1}
\begin{split}
&d\o x_t=[a_t \o x_t+{\o a}_t\o \mu_t -\frac{b_t^2}{n_t} \o y_t+\beta_t]dt
\\
&d\o y_t=-[a_t \o y_t+m_t(m_t \o x_t +{\o m}_t\o \mu_t)]dt,\qquad
\o y_T= q(q \o x_T +{\o q}\o \mu_T),
\end{split}
\end{equation}
has a unique solution. Indeed, 
\eqref{eq:31:12:1} is of the same structure as \eqref{fo:xyfbsde}, with $\sigma=0$ but with the right sign conditions for 
$(b_{t}^2/n_{t})_{0 \leq t \leq T}$, for 
$(m_{t}^2)_{0 \leq t \leq T}$ and for $q^2$, which is enough to repeat the analysis 
from \reff{fo:ansatz:0} to \reff{eq:21:6:2}. 

Similarly, the key point for the solvability of \eqref{eq:29:12:10} is to reformulate it as the adjoint system of a deterministic control problem of linear-quadratic type. Specifically, the minimization problem of the cost function 
\begin{equation*}
\o J(\u{\o \alpha}) = \frac{1}{2} e_{T}q(q+\o q) \xi_{T}^2
+
\int_{0}^T \frac{1}{2} \bigl[ e_{t}n_{t} \o \alpha_{t}^2 + e_{t}m_{t}(m_{t} + \o m_{t}) \xi_{t}^2 \bigr] dt,
\end{equation*}
subject to  
\begin{equation*}
d \xi_{t}= [ (a_{t} + \o a_{t}) \xi_{t} + b_{t} \alpha_{t} + \beta_{t}] dt
\qquad \textrm{and} \qquad
e_{t} = \exp\biggl(-\int^t_{0} \o a_{s} ds\biggr),
\end{equation*}
admits as adjoint system:
\begin{equation}
\label{eq:29:12:11}
\begin{split}
&d \o \xi_{t} =  [(a_{t} + \o a_{t}) \o \xi_{t} - \frac{b_{t}^2}{ e_{t} n_{t}} \o \zeta_{t} + \beta_{t}] dt
\\
&d \o \zeta_{t} = - [e_{t} m_{t} (m_{t}  + \o m_{t})\o \xi_{t} + ( a_{t}+ \o a_{t}) \o \zeta_{t}] dt,
\quad \o \zeta_{T} = e_{T} q ( q+ \o q) \o \xi_{T}.
\end{split}
\end{equation}
Clearly, $(\u {\o \xi},\u {\o \zeta})$ is a solution of \eqref{eq:29:12:11} if and only if 
$(\u {\o \xi},\u e^{-1} \u {\o \zeta})$ is a solution of \eqref{eq:29:12:10}. Unique solvability of \eqref{eq:29:12:11}
is known to hold in short time. In the case when $q(q+\o q) \geq 0$ and $\inf_{t \in [0,T]} [m_{t}(m_{t}+\o m_{t})] \geq 0$,
the optimization problem has a unique optimal path over any arbitrarily prescribed duration $T$ and the associated adjoint system \eqref{eq:29:12:11} is uniquely solvable in finite time as well. In both cases, the fixed point for $\o {\u \mu}$ does exist and is unique:

\begin{theorem}
\label{thm:30:12:1}
On the top of the hypotheses already stated in Subsection \ref{subse:solvNPG},
assume that $q(q+\o q) \geq 0$ and $\inf_{t \in [0,T]} [m_{t}(m_{t}+\o m_{t})] \geq 0$, then, given an initial condition 
for $\u x$, there exists a unique continuous deterministic function $[0,T] \ni t \hookrightarrow \o \mu_{t}$ such that system
\eqref{fo:adjointmfg}--\eqref{fo:fmfg} admits a unique solution $(\underline{x},\underline{y})$ satisfying 
$\o \mu_{t} = \EE(x_{t})$ for any $t \in [0,T]$.
\end{theorem}



\begin{remark}
\label{re:ricatti}
We claimed that limiting ourselves to the one dimensional case $d'=1$ was not restrictive, but it is only fair to mention that, in the higher dimensional framework, the solution of matrix valued 
Riccati's equations is more involved and requires extra assumptions and that the change of variable $(\u {\o \xi},\u {\o \zeta})
\hookrightarrow (\u {\o \xi},\u e^{-1} \u {\o \zeta})$ used in \eqref{eq:29:12:11} 
to reformulate \eqref{eq:29:12:10} cannot be used in such a straightforward way.
\end{remark}

We just explained that a fixed point for $\u{\o \mu}$ could be constructed by investigating the unique solvability of a linear forward-backward ordinary differential equation. Going back to the general discussion \eqref{fo:xyfbsde}--\eqref{fo:mfgfinal}, the solvability of a linear forward-backward equation of deterministic or stochastic type (as in \eqref{fo:xyfbsde}) reduces to the solvability of a Riccati equation of the same type as the first equation in \eqref{fo:mfgfinal}. In practice, once $(\eta_t)_{0 \leq t \leq T}$
in \eqref{fo:mfgfinal} is computed, the point is to plug its value in the third equation to determine $(z_t)_{0 \leq t \leq T}$, and in the second equation, which can then be solved by:
\begin{equation}
\label{eq:30:12:12}
\chi_t= \fr e^{\int_t^T[\fa_u+\fb_u\eta_u]du}-\int_t^T[\fd_s-\fc_s\eta_s]e^{\int_t^s [\fa_u+\fb_u\eta_u]du}ds.
\end{equation}
(When the Riccati equation is well-posed, its solution does not blow up and all the terms above are integrable.)
Now that the deterministic functions $(\eta_t)_{0 \leq t \leq T}$ and $(\chi_t)_{0 \leq t \leq T}$ are computed, we rewrite the forward stochastic differential equation for the dynamics of the state using the ansatz
\eqref{fo:ansatz:0}:
$$
 dx_t=[(\fa_t+\fb_t\eta_t)x_t+\fb_t\chi_t+\fc_t]dt+\sigma dW_t,\qquad x_0=x^{(0)}.
$$
Such a stochastic differential equation is solved explicitly
\begin{equation}
\label{fo:xt}
x_t=x^{(0)} e^{\int_0^t(\fa_u+\fb_u\eta_u)du} +\int_0^t(\fb_s\chi_s+\fc_s)e^{\int_s^t(\fa_u+\fb_u\eta_u)du}ds +\sigma\int_0^te^{\int_s^t
(\fa_u+\fb_u\eta_u)du}dW_s.
\end{equation}

\subsubsection*{A Simple  Example}
For the sake of illustration, we consider the frequently used example where the drift $b$ reduces to the control, namely  $b(t,x,\mu,\alpha)=\alpha$, so that $a_t={\o a}_t=\beta_t=0$, $b_t=1$ and the state equation reads
$$
dx_t=\alpha_tdt+\sigma dW_t, \qquad t \in [0,T]; \qquad x_{0}=x^{(0)}.
$$
We also assume that the running cost is simply the square of the control, i.e. 
$f(t,x,\mu,\alpha)=\alpha^2/2$ and 
$n_t=1$ and $m_t={\o m}_t=0$.
Using the notations and the results above, we see that the FBSDE of the MFG approach has the simple form
\begin{equation}
\label{eq:23:6:1}
\begin{cases}
&dx_t=-y_tdt+\sigma dW_t,
\\
&dy_t=z_tdW_t,
\end{cases} \qquad t \in [0,T]; \qquad  x_{0} = x^{(0)}, \ y_{T} =\fq x_{T} + \fr,
\end{equation}
which we solve by postulating $y_t=\eta_t x_t+\chi_t$, and solving for the two deterministic functions $\eta$ an $\chi$. We find:
\begin{equation*}
\eta_t=\frac{\fq}{1+\fq(T-t)},\qquad\qquad \chi_t=\frac{\fr}{1+\fq(T-t)},
\end{equation*}
(keep in mind that $\fq \geq 0$ so that the functions above are well-defined)
and plugging these expressions into \reff{fo:xt} we get
\begin{equation}
\label{fo:simplext}
x_t=x^{(0)}\frac{1+\fq(T-t)}{1+\fq T}-\frac{\fr t}{1+\fq T}+\sigma[1+\fq (T-t)]\int_0^t\frac{dW_s}{1+\fq (T-s)}.
\end{equation}
Notice further that the optimal control $\alpha_t$ and the adjoint process $y_t$ satisfy
$$
-\alpha_t=y_t=\frac{\fq}{1+\fq (T-t)}x_t+\frac{\fr}{1+\fq(T-t)}
$$
and that the only quantity depending upon the fixed mean function $t\hookrightarrow \o \mu_t$ is the constant $\fr=q{\o q} \, \o \mu_T$, which depends only upon the mean state at the end of the time interval. Recalling that 
$\fq  = q^2$, this makes the search for a fixed point very simple and one easily check that if
\begin{equation}
\label{eq:30:12:10}
\o \mu_T=\frac{x^{(0)}}{1+q(q+{\o q})T}
\end{equation}
then the mean at time $T$ of the random variable $x_T$ given by \reff{fo:simplext} is $\o \mu_T$. 

\begin{remark}
\label{rem:30:12:1}
From \eqref{eq:30:12:10}, we deduce that a fixed point does exist (and in such a case is unique) if 
$1+q(q+ \o q) T >0$. In the case when $q(q+ \o q) \geq 0$, this is always true, as announced in Theorem \ref{thm:30:12:1}. If 
$q(q+\o q) 
<0$, the condition is satisfied if $T$ is small enough only.
\end{remark}
\noindent

\section{Control of Mean Field LQ McKean Vlasov Dynamics}
\label{se:LQMKV}
As before, we restrict ourselves to the one-dimensional case. The problem of the optimal control of the Linear Quadratic McKean-Vlasov dynamics consists in the minimization of the functional
\begin{equation}
\label{fo:MKVobjective}
J(\alpha)=\EE \left[ \int_0^T \bigl[\frac12 \bigl(m_t x_t+\o m_t\EE\{x_t\} \bigr)^2+\frac12 n_t \alpha_t^2\bigr]dt
+\frac 12 \bigl(q x_T+\o q\EE\{x_T\} \bigr)^2\right]
\end{equation}
over all the admissible control processes $\alpha=(\alpha_t)_{0\le t\le T}$ under the constraint
\begin{equation}
\label{fo:McKeanVlasov}
d x_t=\bigl[a_t x_t+\o a_t\EE\{x_t\} + b_t\alpha_t+\beta_t \bigr]dt+\sigma dW_t, \qquad t \in [0,T]; \qquad x_{0}= x^{(0)}.
\end{equation}
As in the case of the MFG approach, we assume that all the coefficients $a_t$, $\o a_t$, $b_t$, $\beta_t$, $ m_t$, $\o m_t$, $ n_t$ are deterministic continuous functions of $t\in[0,T]$, and that $q$  and $\o q$ are deterministic constants. We also assume the function $n_{t}$ to be positive. So in the notation of the Pontryagin principle introduced earlier (recall Theorem \ref{th:pontryaginmkv}), we have
\begin{equation*}
\begin{split}
&\psi(x)= \gamma(x)= \zeta(x)=x, \quad 
b(t,x,x',\alpha)=a_t x +\o a_tx' + b_t \alpha +\beta_t, 
\\
&f(t,x,x',\alpha)=\frac12( m_t x +\o m_tx')^2 +\frac12 n_t \alpha^2, \quad
g(x,x')=\frac12(q x +\o qx')^2,
\end{split}
\end{equation*}
for $(t,x,x') \in [0,T] \times \RR \times \RR$, so that the adjoint equation becomes
\begin{equation}
\label{fo:adjointmkv}
\begin{split}
&dy_t=- \bigl[a_ty_t+ m_t( m_t x_t +\o m_t\EE\{x_t\}) \bigr]dt
- \bigl[\o a_t \EE\{y_t\} +\o m_t(m_t  +\o m_t) \EE\{x_t\} \bigr]dt +z_tdW_t,
\\
&y_T=q \bigl(q x_T +\o q\EE\{x_T\}\bigr) +  \o q (q+ \o q )\EE\{x_{T}\}.
\end{split}
\end{equation}
Notice that because of the special form of the coefficients, this BSDE does not involve the control $\alpha_t$ explicitly. The control appears only indirectly, through the state $x_t$ and its mean $\EE\{x_t\}$. 
Notice also that the sign condition on the derivatives $\partial_{x'} f$ and $\partial_{x'} g$ in (A2) in 
Theorem \ref{th:pontryaginmkv} 
is not satisfied since these two are linear functions. A careful look at the proof of 
Theorem \ref{th:pontryaginmkv} shows that the sign condition is actually useless when the functions $\gamma$ and $\zeta$ are linear, as it is the case here. In order to apply the Pontryagin minimum principle, we guess the candidate for the optimal control by minimizing the Hamiltonian:
\begin{equation*}
\begin{split}
\hat \alpha(t,x,x',y) &=\text{arg}\min_{\alpha\in \RR}H(t,x,x',y,\alpha),
\\
&= \text{arg}\min_{\alpha \in \RR}\biggl\{ y \bigl[a_t x +\o a_t x' +b_t \alpha +\beta_t \bigr]+\frac12( m_t x +\o m_tx')^2 +\frac12 n_t \alpha^2\biggr\}.
\end{split} 
\end{equation*}
Notice that, here, there is one and only one variable $x'$ in the Hamiltonian 
and not two variables $x_{1}'$ and $x_{2}'$ as in \eqref{fo:hamiltonianmkv}. The reason is that the functions $\psi$, $\gamma$
and $\zeta$ coincide. The first order condition gives:
$$
y b_t+ n_t \hat{\alpha}(t,x,x',y)=0, \qquad (t,x,x',y) \in [0,T] \times \RR^3,
$$
so that we choose
\begin{equation}
\label{fo:openalpha}
\hat{\alpha}_t=-\frac{ b_t}{ n_t}y_t, \qquad t \in [0,T],
\end{equation}
as candidate for the optimal control. With this choice, the dynamics of the optimal state become:
\begin{equation}
\label{fo:fmkv}
dx_t= \bigl[a_t x_t+\o a_t\EE\{x_t\} -\frac{ b_t^2}{ n_t}y_t+\beta_t \bigr]dt+\sigma dW_t, \qquad t \in [0,T]; 
\qquad x_{0}=x^{(0)},
\end{equation}
and the problem reduces to the proof of the existence of a solution to the following 
FBSDE of the McKean--Vlasov type (and to the analysis of the properties of such a solution):
\begin{equation}
\label{fo:mkvFBSDE}
\begin{cases}
&\displaystyle dx_t=\bigl[a_t x_t+\o a_t\EE\{x_t\} -\frac{ b_t^2}{n_t}y_t+\beta_t\bigr]dt+\sigma dW_t,
\qquad x_0=x^{(0)},
\\
&dy_t =- \bigl[ a_t y_t+ m_t( m_t x_t +\o{ m}_t\EE\{x_t\}) \bigr]dt 
\\
&\hspace{50pt}- \bigl[ \o a_t \EE\{y_t\} + \o m_t( m_t  +\o m_t) \EE\{x_t\}\bigr]dt
+z_tdW_t,
\\
&y_T=q \bigl(q x_T +\o q\EE\{x_T\}\bigr) +  \o q (q+ \o q )\EE\{x_{T}\}.
\end{cases}
\end{equation}
Since these equations are linear, we can solve first for the functions $t\hookrightarrow \EE\{x_t\}$ and $t\hookrightarrow \EE\{y_t\}$.

Indeed, taking expectations on both sides of \reff{fo:mkvFBSDE} and rewriting the resulting system 
--using the notation $\o x_t$ and $\o y_t$ for the expectations $\EE\{x_t\}$ and $\EE\{y_t\}$ respectively--, we obtain:
\begin{equation}
\label{fo:fbmkv}
\begin{cases}
d\o x_t&\displaystyle  =\bigl[(a_t+\o a_t)\o x_t -\frac{ b_t^2}{ n_t}\o y_t+\beta_t\bigr]dt,\qquad \o x_0= x^{(0)},
\\
\displaystyle
d\o y_t&=- \bigl[(a_t+\o a_t)\o y_t+( m_t +\o m_t)^2 \o x_t \bigr]dt,\qquad
\o y_T=(q +\o q)^2\o x_T.
\end{cases}
\end{equation}
This system is of the same type as \eqref{fo:xyfbsde}, with $\sigma =0$ therein, but with the right sign conditions for 
$(\fb_{t})_{0 \leq t \leq T}$, $(\fm_{t})_{0 \leq t \leq T}$ and $\fq$. In particular, we know from the previous analysis 
 (\ref{fo:ansatz:0}--\ref{eq:21:6:2}) that \eqref{fo:fbmkv} is uniquely solvable. The solution has the form $\o y_t=\o \eta_t \o x_t +\o \chi_t$
for two deterministic functions $\o\eta_t$ and $\o\chi_t$ satisfying
\begin{equation}
\label{fo:meanricatti}
\begin{cases}
&\dot{\o \eta}_t=\displaystyle \frac{b_t^2}{ n_t} \o \eta_t^2 -2 (a_t+\o a_t) \o \eta_t -(m_t +\o m_t)^2,\qquad \o \eta_T=(q +\o q)^2
\vspace{2pt}
\\
&\dot{\o \chi}_t\displaystyle -[\frac{b_t^2}{n_t} \o \eta_{t}-(a_t+\o a_t)] \o \chi_t=-\beta_t \o \eta_t,\qquad
\chi_T=0.
\end{cases}
\end{equation}

Once $\o\eta_t$ is computed, we plug its value in the second equation in \eqref{fo:meanricatti}, which can then be solved
as in \eqref{eq:30:12:12}.
The deterministic functions $\o \eta_t$ and $\o \chi_t$ being computed, we solve for $\o x_t$ and $\o y_t$ by plugging our ansatz $\o y_t=\o \eta_t \o x_t +\o \chi_t$ into the first equation of 
\reff{fo:fbmkv}, as done in \eqref{fo:xt}.
\vskip 6pt

We can now replace the expectations $\EE\{x_t\}$ and $\EE\{y_t\}$ appearing in the FBSDE \reff{fo:mkvFBSDE} of the McKean--Vlasov type by the deterministic functions $\o x_t$ and
$\o y_t$ obtained by solving the forward backward system of ODEs \eqref{fo:fbmkv}, and solve the linear FBSDE
\begin{equation}
\label{fo:mkvFBSDE:1}
\begin{cases}
dx_t&=[\fa_t x_t+\fb_t y_t+\fc_t]dt+\sigma dW_t,\qquad x_0=x^{(0)}\\
dy_t&=[\fm_t x_t- \fa_t y_t+ \fd_t]dt+z_tdW_t,\qquad
y_T= \fq x_T +\fr,
\end{cases}
\end{equation}
where for the purpose of this part of the proof, and in order to lighten the notation, we use the same kind of notation we used in the case of the MFG approach by setting:
\begin{equation*}
\begin{split}
&\fa_t=a_t, \;  \fb_t=- b_t^2/ n_t,\; \fc_t=\beta_t+\o a_t\o x_t, 
\\
&\fm_t= - m_{t}^2, \; \fd_t=-\o m_t(2 m_t+\o m_t)\o x_t-\o a_t \o y_t, \;
\fq=q^2,\; \fr=\o q(2q+\o q)\o x_T.
\end{split}
\end{equation*}
Here as well, we are in the framework of equation \reff{fo:xyfbsde}, with the right-sign conditions, so that equation \eqref{fo:mkvFBSDE:1} is uniquely solvable. The pair process $(\EE\{x_{t}\},\EE\{y_{t}\})_{0 \leq t \leq T}$ is then the solution to an ordinary forward-backward system, which is also of the same type as \reff{fo:xyfbsde}: the solution must be the pair $(\o x_{t},\o y_{t})_{0 \leq t \leq T}$ so that \eqref{fo:mkvFBSDE} is uniquely solvable. Clearly, the FBSDE value function is affine, that is 
 \begin{equation}
\label{fo:ansatz}
y_t=\eta_t x_t+\chi_t,\qquad\qquad 0\le t\le T,
\end{equation}
for some $(\eta_{t})_{0 \leq t \leq T}$ and $(\chi_{t})_{0 \leq t \leq T}$ given by a system of the same type 
as \eqref{fo:mfgfinal}. We claim:
\begin{proposition}
\label{prop:30:12:1}
Under the current assumptions, for a given initial condition for the forward equation, the FBSDE
\reff{fo:mkvFBSDE} is uniquely solvable.
\end{proposition}

Once $(\eta_t)_{0 \leq t \leq T}$ and $(\chi_t)_{0 \leq t \leq T}$ in \eqref{fo:ansatz} are determined by solving the corresponding Riccati equation, we rewrite the forward stochastic differential equation for the dynamics of the state:
$$
dx_t=\bigl[(\fa_t+\fb_t\eta_t)x_t+\fb_t\chi_t+\fc_t\bigr]dt+\sigma dW_t, \qquad t \in [0,T]; \qquad x_0=x^{(0)},
$$
which is solved explicitly, as in \eqref{fo:xt}.

\vskip 2pt
Notice that \reff{fo:xt} shows that the optimally controlled state is still Gaussian despite the nonlinearity due to the McKean-Vlasov nature of the dynamics. While the expectation $\EE\{x_t\}=\o x_t$ was already computed, expression \reff{fo:xt} can be used to compute the variance of $x_t$ as a function of time. Because of the linearity of the ansatz and the fact that $\eta_t$ and $\chi_t$ are deterministic,
the adjoint process $(y_t)_{0\le t\le T}$ is also Gaussian.

\begin{remark}
Using again the form of the ansatz, we see that the optimal control $\alpha_t$ which was originally identified as an open loop control in \reff{fo:openalpha}, is in fact in closed loop feedback form since it can be rewritten as
\begin{equation}
\label{fo:closealpha}
\alpha_t=-\frac{b_t}{n_t}\eta_t x_t-\frac{b_t}{n_t}\chi_t
\end{equation}
via the feedback function $\phi(t,\xi)=-b_t(\eta_t\xi+\chi_t)/ n_t$ which incidently shows that the optimal control is also a Gaussian process.
\end{remark}

\begin{remark}
The reader might notice that, within the linear-quadratic framework, the conditions for the unique solvability of the adjoint equations are not the same in the MFG approach and in the control of MKV dynamics. On the one hand, optimization over controlled MKV dynamics reads as an optimization problem of purely (strictly) convex nature, for which existence and uniqueness of an optimal state is expected. On the other hand, the optimal states in the MFG approach appear as the fixed points of a matching problem. Without any additional assumptions, there is no reason why this matching problem should derive from a convex potential even if the original coefficients of the game are linear-quadratic. In Theorem \ref{thm:30:12:1}, the sign conditions on 
$q(q+\o q) \geq 0$ and $\inf_{t \in [0,T]} [m_{t}(m_{t}+\o m_{t})] \geq 0$ are precisely designed to make the matching problem 
\eqref{eq:31:12:1} be (strictly) convex.
\end{remark}

\subsubsection*{Simple  Example}
We consider the same example as in the MFG approach where $b(t,x,\mu,\alpha)=\alpha$, so that $a_t={\o a}_t=\beta_t=0$, $b_t=1$ so that the state equation reads
$$
dx_t=\alpha_tdt+\sigma dW_t, \qquad t \in [0,T]; \qquad x_{0}=x^{(0)},
$$
and $f(t,x,\mu,\alpha)=\alpha^2/2$, so that $n_t=1$ and $m_t={\o m}_t=0$.
Using the notation and the results above, we see that the system 
\eqref{fo:mkvFBSDE} has the same form as \eqref{eq:23:6:1}, but $(\fq,\fr)$ in \eqref{eq:23:6:1} is now given by
$\fq = q^2$ and $\fr =  \o q (2q + \o q) \EE\{x_{T}\}$. Postulating the relationship $y_t=\eta_t x_t+\chi_t$ and solving for the two deterministic functions $\eta$ and $\chi$, we find the same expression as in \eqref{fo:simplext}:
\begin{equation*}
x_t=x^{(0)}\frac{1+\fq(T-t)}{1+\fq T}-\frac{\fr t}{1+\fq T}+\sigma[1+\fq (T-t)]\int_0^t\frac{dW_s}{1+\fq (T-s)}.
\end{equation*}
This makes the computation of $\EE\{x_{T}\}$ very simple. We find 
\begin{equation*}
\EE\{x_T\}=\frac{x^{(0)}}{1+(q+{\o q})^2T},
\end{equation*}
which always makes sense and which is different from \eqref{eq:30:12:10}.

\section{Further Examples}
\label{se:examples}
In this section we provide simple examples which do not fit in the linear quadratic mold of the previous sections, but still allow for an explicit comparison of the solutions of the two problems. Surprisingly, we shall see that these solutions do coincide in some special cases. 

\subsection{Some Simple Particular Cases}

As in previous sections, we discuss the case when $b(t,x,\mu,\alpha)=\alpha$ and 
$f(t,x,\mu,\alpha)=\alpha^2/2$. When the terminal cost is of quadratic type, the previous analysis shows that, in both cases, the optimal control is simply the negative of the adjoint process, the difference between the two cases lying in the shape of the terminal condition of the associated FBSDEs. We discuss now several related forms of the function $g$ for which we can still conclude.

\bigskip
\noindent
{\bf The case $g(x,\mu)= r x \o {\mu}$, $r \in \RR^*$.} Here we only compare the end-of-the-period means. We know from 
Section \ref{se:LQMFG} that the solution(s) to the MFG problem are given by the FBSDE \eqref{eq:23:6:1} with the terminal condition 
$\fq=0$ and $\fr =  r \o{\mu}_{T}$, so that, from \reff{fo:simplext}, the fixed point condition reads
\begin{equation}
\label{fo:mean1}
\o \mu_{T} = x^{(0)} -  r \o{\mu}_{T} T, \qquad \textrm{i.e.} \quad 
(1+rT) \o \mu_T=x^{(0)},
\end{equation}
which says that there exists a unique fixed point when $1+rT \not = 0$. If $1+rT=0$, then there is no fixed point if $x^{(0)} \not = 0$ and there are infinitely 
many if $x^{(0)}=0$. Given $\o \mu_{T}$ as in \eqref{fo:mean1}, we know from the analysis in Section \ref{se:LQMFG} that \reff{eq:23:6:1} has a unique solution $(x_{t},y_{t},z_{t})_{0 \leq t \leq T}$. By \reff{fo:simplext}, it must hold 
$\EE\{x_{T}\}=\o \mu_{T}$, so that $(x_{t},y_{t},z_{t})_{0 \leq t \leq T}$ is the unique equilibrium of the MFG problem. 

When dealing with the control of MKV dynamics with the same choice of coefficients, the story is rather different. Indeed, the terminal condition $g$ is not convex in the variables $(x,x')$, so that Theorem \ref{th:pontryaginmkv} does not apply anymore. From 
\cite{AndersonDjehiche}, we know that the FBSDE \eqref{fo:adjointmkv:0} just provides a necessary condition for the optimal states of the optimization problem. We thus investigate the solutions to \eqref{fo:adjointmkv:0} as possible candidates for minimizing the related cost. As above, \eqref{fo:adjointmkv:0} reduces to \eqref{eq:23:6:1}, with the terminal condition 
$\fq = 0$ and $\fr = 2 r \EE\{x_{T}\}$. By \eqref{fo:simplext}, the necessary condition for 
$\EE\{x_{T}\}$ reads
\begin{equation*}
\EE\{x_{T}\} = x^{(0)} - 2 r \EE\{x_{T}\} T, \qquad \textrm{i.e.} \quad 
(1+2rT) \o \mu_T= x^{(0)},
\end{equation*}
which says that there exists a unique fixed point when $1+2rT \not =0$. When $1+2r T = 0$, there is no fixed point unless
$x^{(0)}=0$. Moreover, the optimal states, when they do exist, must differ from 
the optimal states of the MFG problem (unless
$x^{(0)}=0$).
\bigskip

\noindent
{\bf The case $g(x,\mu)= r x \o \mu^2$, $r \in \RR^*$.} We argue as above. In the MFG framework, 
the fixed point condition also derives from the FBSDE \eqref{eq:23:6:1}, with the terminal condition 
$\fq=0$ and $\fr =  r \o{\mu}_{T}^2$, and from \reff{fo:simplext}. It reads
\begin{equation}
\label{meanmfg}
\o \mu_{T} = x^{(0)} - r T \o{\mu}_{T}^2, \qquad \textrm{i.e.} \quad
r T \o{\mu}_{T}^2 + \o \mu_{T} - x^{(0)}  =0.
\end{equation}
By the same argument as in the previous case, given a solution to \eqref{meanmfg}, 
there exists a unique solution to \eqref{eq:23:6:1}. In particular, the number of solutions to the MFG problem matches the number of solutions to \eqref{meanmfg} exactly.

For the MKV control problem, the FBSDE \eqref{eq:23:6:1}, with the terminal condition 
$\fq=0$ and $\fr =  3 r \EE\{x_{T}\}^2$, just provides a necessary condition for the existence of optimal states. 
By \reff{fo:simplext}, $\EE\{x_{T}\}$ must solve:
\begin{equation}
\label{meanmkv}
\EE \{x_{T}\} = x^{(0)} - 3 r T \EE\{x_{T}\}^2, \qquad \textrm{i.e.} \quad
3 r T \EE \{x_{T}\}^2 + \EE\{x_{T}\} - x^{(0)}  =0.
\end{equation}
Existence of a solution and uniqueness to the second-order-equations \reff{meanmfg} and \reff{meanmkv} depend upon the values of the parameters $r$, $T$, and $x^{(0)}$. Let us have a quick look at these issues  in more detail. Let $S_{\textrm{MFG}}=\{(T,x) \in \RR_+ \times \RR; \; 1+4 r Tx^{(0)}\geq 0\}$ and $S_{\textrm{MKV}}
=\{(T,x^{(0)}) \in \RR_+ \times \RR; \; 1+12 r Tx^{(0)}\geq 0\}$. These represent the sets where the
second-order-equations are solvable. It is easy to see that $S_{\textrm{MFG}}\subset S_{\textrm{MKV}}$ when $r x^{(0)}>0$, with the converse inclusion when $r x^{(0)}<0$. Moreover, there exists a continuum of triplets $(r,T,x^{(0)})$ for which there is at most one solution  to the MKV control problem while there are exactly two solutions  to the MFG problem. In particular, if $r=- 1$, $T=1$, $x^{(0)}=1/12$, then we have $(1,1/6) \in S_{\textrm{MKV}}$ and $(1,(3\pm \sqrt{6})/6) \in S_{\textrm{MFG}}$.

\noindent {\bf General Linear Terminal Cost.} The two previous examples can be easily generalized to the case when the terminal cost is of the form $g(x,\mu)=x \gamma(\bar{\mu})$. In the MFG approach, the fixed point equation for the mean has the form (compare with 
\eqref{fo:mean1}):
\begin{equation*}
\o \mu_T=x^{(0)} - T \gamma(\o \mu_T).
\end{equation*}
For the MKV control problem, it reads
\begin{equation*}
\EE\{x_T\}=x^{(0)} - T \bigl[ \gamma'\bigl(\EE\{x_T\}\bigr)
\EE\{x_T\}  + \gamma\bigl(\EE\{x_T\} \bigr) \bigr].
\end{equation*}

\noindent {\bf Quadratic Terminal Cost.} When $g(x,\o \mu)=x^2  \gamma(\o \mu)$, we deduce, by choosing $\fq = 2 \gamma(\o \mu)$  and $\fr=0$ in \eqref{eq:23:6:1}, that, in the MFG approach,
the fixed point equation for the mean has the form 
\begin{equation*}
\o \mu_T \bigl( 1+ 2 T \gamma(\o \mu_{T}) \bigr) =x^{(0)}.
\end{equation*}
For the MVK control problem, the analysis of the adjoint equation cannot be tackled by investigating the mean of $x_{T}$ first since the terminal condition depends upon the expectation of the square of $x_{T}$.

\subsection{Additive Running Cost} Here we tackle two examples when $b(t,x,\mu,\alpha)=\alpha$ and $f \not = 0$. In both cases, 
the terminal cost $g$ is set equal to $0$. 

\bigskip
\noindent
{\bf The case $f(t,x,\mu,\alpha)= \alpha^2/2 + x \o \mu$.} In the MFG approach, we can write the equivalent of the averaged forward-backward differential equation
\eqref{eq:31:12:1}. We get as a fixed point condition for the mean $(\o \mu_{t})_{0 \leq t \leq T}$:
\begin{equation*}
\begin{cases}
&d \overline{\mu}_{t} = - \o p_{t} dt, 
\\
&d\o p_{t} = - \o \mu_{t} dt, 
\end{cases}
\qquad t \in [0,T]; \quad \o \mu_{0} = x^{(0)}, \quad \o p_{T}=0.
\end{equation*}
This reads as the second-order ODE $\ddot{\o \mu}_{t} = \o \mu_{t} = 0$, with
$\o \mu_{0} = x^{(0)}$ and  $\dot{\o \mu}_{T} = 0$,
the solution of which is given by 
\begin{equation*}
\o \mu_t=x^{(0)}\frac{e^{(T-t)}+e^{-(T-t)}}{e^{T}+e^{-T}}, \quad t \in [0,T].
\end{equation*}
As in the simple case investigated in the previous subsection, once $(\o \mu_{t})_{0 \leq t \leq T}$ is computed, the solution to 
the corresponding FBSDE \eqref{eq:16:12:2} exists and is unique, since the FBSDE is then of the same form 
as \eqref{fo:xyfbsde}.

For the MKV problem, convexity does not hold, so that the stochastic minimum principle provides a necessary condition only. Here it has the form:
\begin{equation*}
\begin{cases}
&d \EE\{x_{t}\} = - \EE\{ y_{t}\} dt, 
\\
&d\EE \{y_{t}\} = - 2\EE \{x_{t}\} dt, 
\end{cases}
\qquad t \in [0,T]; \quad  \EE\{x_{0}\} = x^{(0)}, \quad \EE\{\o y_{T}\}=0.
\end{equation*}
Again, the above system reduces to a second-order ODE, the solution of which is:
\begin{equation*}
\EE \{x_t\}=x^{(0)}  \frac{e^{\sqrt{2}(T-t)}+e^{-\sqrt{2}(T-t)}}{e^{\sqrt{2}T}+e^{-\sqrt{2}T}}, \quad t \in [0,T].
\end{equation*}

\noindent
{\bf The linear-quadratic case with zero terminal cost.} We finally go back to the LQ framework. As above, we choose $b(t,x,\mu,\alpha)=\alpha$ and $g=0$. For the running cost, we consider the classical case where 
$f(t,x,\mu,\alpha)= \alpha^2/2 + (x-\o \mu)^2/2$. Here it is easy to see that, in both cases, the solution to the adjoint backward equation has zero mean, so that the expectation of any optimal states must be constant.

\subsection{A Simple Model for Emissions Regulation}
\label{se:co2}
For the sake of motivation, we present a toy model of the simplest form of Green House Gas (GHG) emissions regulation. Our interest in this model
is the fact that the mean field interaction between the players appears naturally in the terminal cost function.

A set $\{1,\dots,N\}$ of firms compete in an economy where green house gas emissions are regulated over a period $[0,T]$, for some $T>0$. For each firm $i$, we denote by $X^i_t$ the perceived value at time $t$ of what its own 
cumulative emissions will be at maturity $T$, and we assume that its dynamics satisfy
\begin{equation}
\label{fo:edynamics}
 dX^i_t=(b^i_t-\alpha^i_t)dt+\sigma^i_tdW_t, \qquad t \in [0,T]; \quad X^i_0= x^{(0)},
 \end{equation}
where assumptions on the individual emission rates $b^i_t$ and the volatilities $\sigma^i_t$ can be chosen as in \cite{CDET} for the sake of definiteness. The process  $(\alpha^i_t)_{0\le t\le T}$ represents the abatement rate of firm $i$, and it can be viewed as the control the firm can exert on its emissions output. In the case $\alpha^i_t\equiv 0$, the process  $(X^i_t)_{0\le t\le T}$ gives the perceived cumulative emissions of firm $i$ in the absence of regulation, situation which is called \emph{Business as Usual}, BAU in short. At the start $t=0$ of the regulation implementation period, each firm $i$ is allocated a number $\Lambda_i$ of permits (also called certificates of emission, or allowances). The cap
$\Lambda^{(N)} = \sum_{i=1}^N \Lambda_i$
can be viewed as the emissions target set by the regulator for the period $[0,T]$. If at the end of the regulation period $[0,T]$, the aggregate emissions in the economy exceed the cap, i.e. if $
\sum_{i=1}^N X^i_T > \Lambda^{(N)}$,
then each firm $i$ has to offset its emissions $X^i_T$ (expressed in CO${}_2$ ton equivalent) by redeeming one permit per ton, or by paying a penalty $\lambda$ for each ton not covered by a certificate. In other words, firm $i$ has to pay
\begin{equation}
\label{fo:cap}
\lambda_i=\lambda (X^i_T-\Lambda_i)^+\bone_{(\Lambda^{(N)} ,\infty)} \biggl(\sum_{j=1}^NX^j_T \biggr)
\end{equation}
where we use the notation $x^+=\max(x,0)$ for the positive part of the real number $x$, and $\bone_A$ for the indicator function of the set $A$.
The penalty $\lambda$ is currently equal to $100$ euros in the European Union Emissions Trading Scheme (EU ETS). 

\begin{remark} In the case of cap-and-trade regulations, a market on which emissions certificates can be bought and sold is also created as part of the regulation.
In order to avoid penalties, firms whose productions are likely to lead to cumulative emissions in excess of their initial allocations $\Lambda_i$ may engage
in buying allowances from firms which expect to meet demand with less
emissions than their own initial allocations. For the sake of simplicity, the present discussion concentrates on the \emph{cap} part, and disregards the \emph{trade} part of the regulation.
\end{remark}

In search for an optimal behavior, each firm needs to solve the following optimization problem.
If we assume that the abatement costs for firm $i$ are given by a function $c^i:\RR\to\RR$
which  is $C^1$ and strictly convex, is normalized to satisfy $c^i(0)=\min c^i$ , and satisfies Inada-like conditions, ($c^i(x)=\beta|x|^{1+\alpha}$  for some $\beta>0$ and $\alpha>0$ is an example of such a cost function), and if its abatement strategy is $\alpha^i$, its terminal wealth is given by
$$
 W^i_T=w^i -\int_0^Tc^i(\alpha^i_t)dt-\lambda(X^i_T-\Lambda_i)^+\bone_{(\Lambda^{(N)},\infty)}
 \biggl(\sum_{j=1}^NX^j_T \biggr),
 $$
 where $w^i$ stands for the initial wealth of the firm $i$.
 Recall that, in our simple model, firms do not trade allowances.
Now if we view the perceived emissions $X^i_t$ as the private state of firm $i$, and if we assume that each firm tries to maximize its expected terminal wealth, or equivalently minimize the objective function
\begin{equation}
\label{fo:ecost}
J^i(\u \alpha)=\EE\biggl[\int_0^Tc^i(\alpha^i_t)dt + \lambda (X^i_T-\Lambda_i)^+\bone_{(\Lambda^{(N)},\infty)}\biggl(\sum_{j=1}^NX^j_T\biggr)\biggr],
\end{equation}
where we set $\alpha_t=(\alpha^1_t,\cdots,\alpha^N_t)$, then we have formulated our emission regulation model as
a stochastic differential game.
Notice that, in this particular example, the equations \reff{fo:edynamics} giving the dynamics of the private states $X^i_t$ are decoupled. However, equations \reff{fo:ecost} giving the expected costs to minimize are coupled in a very specific way, namely through the average of the values of all the private states. 
Recasting the model in the framework of Section \ref{se:games}, we see that the drift and running costs of the private states $X^i_t$ are of the form:
$$
b(t,x,\mu,\alpha)=-\alpha,\qquad \text{and}\qquad f(t,x,\mu,\alpha)=c(\alpha).
$$
For the sake of simplicity we assume that the BAU drift $b^i_t$ is zero, even if, from a modeling point of view, this implies that the processes $X^i$ can take negative values. We refer the reader to \cite{CDET} for a specific discussion of the positivity of the emissions in similar models. The terminal cost is given by the function $g$ defined on $\RR\times \cP_1(\RR)$ by
$$
g(x,\mu)=\lambda (x-\Lambda)^+ \bone_{\{\o\mu>\Lambda\}}
$$
where $\Lambda=\Lambda^{(N)}/N$ stands for the cap per firm, which is assumed to be independent of $N$. The quantity $\Lambda$ is the relevant form of the cap as it makes sense in the limit $N\to\infty$ of a large number of firms treated similarly by the regulator.

\subsubsection*{The Mean Field Game Approximation}

If we assume quadratic abatement costs $c(\alpha)=\alpha^2/2$, then the approximate Nash equilibriums provided by the MFG approach are given by the solutions $(x_{t},y_{t},z_{t})_{0 \leq t \leq T}$ 
of the FBSDE \eqref{eq:23:6:1}, with $0$ as initial condition and terminal condition 
\begin{equation*}
y_{T} = \lambda {\mathbf 1}_{\{x_{T} > \Lambda \}} {\mathbf 1}_{\{\o \mu_{T} > \Lambda \}},
\end{equation*}
satisfying the matching problem $\EE\{x_{T}\} = \o \mu_{T}$, that is
\begin{equation}
\label{eq:30:6:1}
\EE\{x_{T}\} = x^{(0)} -  T \lambda {\mathbb P}\{x_{T}>\Lambda\}
{\mathbf 1}_{\{\EE\{x_{T}\} > \Lambda \}}. 
\end{equation}
We refer the reader to \cite{CDET} for the
precise derivation of the terminal condition despite the singularity at $\Lambda$. Because of the nonlinearity of the terminal condition, the matching problem is much more involved than in the linear-quadratic setting. Precisely, the simultaneity of the nonlinearity of the terminal condition and of the strong coupling between the forward and backward SDEs makes the fixed point equation rather intricate. In the specific situation when the running cost is quadratic, as it is here assumed to be, the forward and backward equations can be decoupled by a Hopf-Cole transformation which suggests to work with an exponential of the solution of the Hamilton-Jacobi-Bellman equation. This transformation also applies to the FBSDE by differentiating the value function of the control problem. Remember that the value function of the FBSDE coincides with the derivative of the value function of the control problem. Therefore it is simpler to investigate the Hamilton-Jacobi-Bellman equation directly. In the present set-up, the  forward-backward system of nonlinear PDEs describing the MFG problem reads (compare with (\ref{hjb}--\ref{fo:kolmo})):
\begin{eqnarray}
\label{mfg}
\mbox{ \footnotesize (Kolmogorov)} \;\; &&\partial_t m -\frac{\sigma^2}{2} \partial_{xx}^2 m + \partial_{x}( - \partial _{x} v m)=0, \; \; \;  \; m(0,\cdot)=\delta_{x^{(0)}},  \nonumber \\ 
\mbox{ \footnotesize (HJB)} \;\; &&\partial_t v + \frac{\sigma^2}{2} \partial_{xx}^2 v- \frac{1}{2}(\partial_{x}v )^2=0, \; \; \; \; \; \; \; \; \; v(T,\cdot)=
\lambda ( \cdot - \Lambda)^+ {\mathbf 1}_{\{\bar{m}_{T}>\Lambda\}}, 
\nonumber
\end{eqnarray}
where the function $m(t,\,\cdot\,)$ is understood as the density of the law of $X_{t}$ in \eqref{eq:16:12:2} (or equivalently as the density of the fixed point $(\mu_{t})_{0 \leq t \leq T}$ of the matching problem), so that:
$$
\overline m_t= \o \mu_{t} = \int_{\RR} x m(t,x)\,dx, \quad t \in [0,T],
$$
and $v$ as the HJB value function of the stochastic control problem when the family $(\mu_t)_{0\le t\le T}$ is frozen. 
These PDEs are simple enough to be solved by a Hopf-Cole transformation:
\begin{equation*}
v=-\sigma^2 \log{u} \leftrightarrow u = \exp\left({-\frac{v}{\sigma^2}}\right) \; ;
\qquad
m= u \psi  \leftrightarrow \psi = \frac{m}{u}. 
\end{equation*}
Then the pair $(u,\psi)$ solves the simple forward-backward system:
\begin{eqnarray}\label{heatsystem}
&&\partial_t \psi -\frac{\sigma^2}{2} \partial_{xx}^2\psi=0,    \; \; \;  \; \; \; \; \; \psi(0,\cdot)= \frac{1}{u(0,x^{(0)})}
\delta_{x^{(0)}},  
\nonumber 
\\ 
&&\partial_t u + \frac{\sigma^2}{2} \partial^2_{xx}u=0, \; \; \;  \; \; \; \; \; u(T,\cdot)=
\exp \bigl[ - \lambda \sigma^{-2} (\cdot - \Lambda)^+ {\mathbf 1}_{\{\bar{m}_{T}>\Lambda\}} \bigr], 
\nonumber
\end{eqnarray}
 consisting of two fundamental heat equations coupled only at the end points of the time interval. 

Now the specific and simple form of the mean-field interaction term $\o m_T$ (depending only upon the global distribution's first moment) leads us to distinguish three cases:
\vspace{2pt}

\textit{The Business As Usual (BAU) Solution.} If we assume that the cap is not reached, then $u(T,y)=1$ for all $y\in \R$, which implies that, for all $(t,y) \in [0,T]\times \RR$,  we have $u(t,y)=1$ and thus $v(T,y)=0$. In this framework, the optimal abatement strategy is to do nothing, i.e. $\hat{\alpha}_{t}=0$. This corresponds to BAU. So, in this case,
$m(t,y) = \varphi_{x^{(0)},\sigma^2t}(y)$, the density of the Gaussian distribution with mean $x^{(0)}$ and variance $\sigma^2 t$.
The fixed point condition for the Nash equilibrium is then satisfied if:
\begin{equation*}
\int_{\RR} y \varphi_{x^{(0)},\sigma^2t}(y) dy \leq \Lambda \qquad \textrm{i.e.} \quad x^{(0)} \leq \Lambda.
\end{equation*}
\vspace{2pt}

\textit{The Abatement Solution.} Suppose now that the cap is exceeded. In this case the terminal condition has the form: $v(T,y)=\lambda (y-\Lambda)^+$. An easy computation shows that 
$$v(t,x)=-\sigma^2 \log\biggl( \int_{\RR} \varphi_{x,\sigma^2(T-t)}(y) \exp \bigl[ - \lambda \sigma^{-2} (y-\Lambda)^+\bigr] dy\biggr), \quad (t,x) \in [0,T] \times \RR.$$
The optimal feedback is given by $\partial_{x} v(t,x)$, which satisfies the viscous Burgers equation:
\begin{equation*}
\partial_t \bigl(\partial_{x}v) + \frac{\sigma^2}{2} \partial_{xx}^2 \bigl( \partial_{x}v \bigr)
- \partial_{x}v \partial_{x}(\partial_{x} v) =0, \; \; \; \; \; \; \; \; \; \partial_{x} v(T,\cdot)=
\lambda {\mathbf 1}_{\{\cdot>\Lambda\}}.
\end{equation*}
By the maximum principle, it holds $0 \leq \partial_{x} v \leq \lambda$. Therefore, when $x^{(0)} > \Lambda + \lambda T$, the optimal path $(x_{t})_{0 \leq t \leq T}$ satisfies at maturity:
\begin{equation}
\label{eq:30:6:2}
x_{T} > \Lambda + \sigma W_{T},
\end{equation}
so that $\bar{m}_{T} > \Lambda$, which guarantees that the fixed point condition for the Nash equilibrium is satisfied.  

\textit{Critical Case.} If $x^{(0)}$ belongs to the interval $(\Lambda, \Lambda+\lambda T/2)$, existence of the fixed point may fail. Indeed \eqref{eq:30:6:1} says that a necessary condition for the existence of a fixed point is $\EE(x_{T})>\Lambda$, so that $x^{(0)}$ must satisfy
\begin{equation*}
x^{(0)} \geq \EE(x_{T}) + \lambda T \PP\bigl\{x_{T} > \Lambda \bigr\}.
\end{equation*}
Following the proof of \eqref{eq:30:6:2}, we see that $x_{T} \geq x^{(0)} - \lambda T + \sigma W_{T}
\geq \Lambda - \lambda T + \sigma W_{T}$. Therefore, 
\begin{equation*}
x^{(0)} \geq \EE(x_{T}) + \lambda T \PP\bigl\{\sigma W_{T} > \lambda T \bigr\}
> \Lambda + \lambda T \PP\bigl\{\sigma W_{T} > \lambda T \bigr\}.
\end{equation*}
When $\sigma$ tends to $+\infty$, the inequality degenerates into 
\begin{equation*}
x^{(0)} 
\geq \Lambda + \frac{\lambda T}{2},
\end{equation*}
so that, when $x^{(0)} < \Lambda + \lambda T/2$,  existence of the fixed point fails for $\sigma$ large enough. 

When $x^{(0)} \geq \Lambda + \lambda T/2$, we must refine the previous argument in order to settle the question of existence of a fixed point. 
The difficulty is to compute $\PP \{ x_{T} > \Lambda \}$. 
Using the expression of $m$ in terms of $\psi$, we see that: 
\begin{equation*}
\PP \bigl\{ x_{T} > \Lambda \bigr\} 
= \frac{\displaystyle \int_{\Lambda}^{+\infty} u(T,y) \varphi_{x^{(0)},\sigma^2 T}(y) dy} 
{\displaystyle \int_{-\infty}^{+\infty} u(T,y) \varphi_{x^{(0)},\sigma^2 T}(y) dy}.
\end{equation*}
Moreover,
\begin{equation*}
\begin{split}
&\int_{\Lambda}^{+\infty} u(T,y) \varphi_{x^{(0)},\sigma^2 T}(y) dy
\\
&= (2 \pi)^{-1/2}\sigma^{-1} T^{-1/2} \int_{0}^{+\infty} 
\exp \bigl[ - \lambda \sigma^{-2} y - \frac{1}{2} \sigma^{-2} T^{-1}
\bigl( y - (x^{(0)}-\Lambda)\bigr)^2 \bigr] dy
\\
&= \exp \bigl[ \lambda \sigma^{-2} \bigl( \lambda T - 2 (x^{(0)}-\Lambda) \bigr)/2 \bigr] \Phi \bigl( \frac{x^{(0)}
- \Lambda - \lambda T}{\sigma T^{1/2}} \bigr),
\end{split}
\end{equation*}
where $\Phi$ stands for the cumulative distribution function of the standard Gaussian distribution, so since
\begin{equation*}
\int_{-\infty}^{\Lambda} u(T,y) 
\varphi_{x^{(0)},\sigma^2 T}(y) = \PP 
\{ \sigma W_{T} \leq    \Lambda -x^{(0)} \} = \Phi 
\bigl( \frac{\Lambda - x^{(0)}}{\sigma T^{1/2}} \bigr),
\end{equation*}
we deduce
\begin{equation}
\label{eq:30:6:7}
\PP \bigl\{ x_{T} > \Lambda \bigr\} 
= \frac{\displaystyle \exp \bigl[ - \lambda \sigma^{-2} \delta \bigr] \Phi \bigl[ 
\sigma^{-1} T^{-1/2} \bigl(- \lambda T/2 + \delta\bigr) \bigr]}{\displaystyle 
\exp \bigl[ - \lambda \sigma^{-2} \delta \bigr] \Phi \bigl[ 
\sigma^{-1} T^{-1/2} \bigl( - \lambda T/2+ \delta\bigr) \bigr] + \Phi
\bigl[ \sigma^{-1} T^{-1/2} \bigl( - \lambda T/2 - \delta\bigr)\bigr]},
\end{equation}
where we have set $\delta =  x^{(0)} - \Lambda - \lambda T/2$. 

When $\delta=0$, $\PP \{ x_{T} > \Lambda \}=1/2$ so that $\EE\{x_{T}\}=\Lambda$. This suggests that 
$x^{(0)} = \Lambda + \lambda T/2$ is a critical initial condition for the existence of a Nash equilibrium. This can be proven rigorously
when $\sigma$ is large. Keep in mind the fact that existence then fails for $x^{(0)} < \Lambda + \lambda T/2$.
Indeed, letting $\sigma$ tend to $+\infty$ in 
\eqref{eq:30:6:7},  we obtain  
\begin{equation*}
\PP \bigl\{ x_{T} > \Lambda \bigr\} \rightarrow \frac{1}{2} 
\end{equation*}
so that, by \eqref{eq:30:6:1},
 $\EE(x_{T}) > \Lambda$ for $\sigma$ large, and the fixed point condition holds. 
\vskip 4pt
An intuitive explanation why the fixed point condition can fail goes as follows. In a reasonable (from the point of view of the choices of $\Lambda$, $\lambda$ and $T$) cap-and-trade scheme, the cap is expected to be reached to incentivize the implementation of abatement strategies. However,  it is often observed that the cap is not reached at the end of the period. From the MFG theory viewpoint, the argument is as follows. An individual firm with negligible impact on the overall emissions can emit whatever it wants without impacting significantly the global emissions. However, as soon as this becomes everybody's strategy, the cap is reached.

\subsubsection*{Search for Cooperative Equilibriums}

From an economical point of view, one might think that the search for an approximate Nash equilibrium is well-adapted to the case when the firms decide of their own strategy without any significant effect on the strategies of the others. Even if this situation sounds reasonable, the reader might wonder about the case when some of the firms obey a similar policy. Part of the answer is then given by 
Subsection 
\ref{subse:COOP}. Indeed, in the extremal situation when the all the firms follow the same general strategy, equilibriums are expected to be of the cooperative type, as described in Subsection 
\ref{subse:COOP}. As the terminal cost function $g$ is non-convex with respect to the expectation $\bar{\mu}$, the stochastic Pontryagin principle then provides necessary conditions only for the optimal strategies. Within this framework, it is then quite remarkable that the forward-backward equation 
\eqref{fo:adjointmkv:0} matches 
\eqref{eq:16:12:2} excatly, provided that the expectation of the emissions at maturity is different from the cap $\Lambda$. (So that both equations admit the same solutions.) The reason why they match is as follows. If the expectation of the emissions at time $T$ 
differs from the cap, any slight perturbation of the overall strategy keeps the final state of the system unchanged: the cap remains either exceeded or respected. Therefore, the search for a local cooperative equilibrium leads to the same limit optimization problem as the search for a local Nash equilibrium, which means that the necessary conditions deriving from the stochastic Pontryagin principle are the same. The case when the mean of the emissions at maturity fits the cap exactly is much more involved and goes far beyond the scope of the paper: as shown in \cite{CDET}, the collateral effect of the cap onto the dynamics of the carbon market is highly singular in the vicinity of the cap when the market has some degeneracy. So is the case here as the dynamics for the mean are completely degenerate. (The previous paragraph shows that the dynamics for the mean are related to the inviscid Burgers equation.)

\section{General Solvability Results}
\label{se:solvability}
We argued that optimal paths to both MFG and MKV control problems appear as realizations of the forward component of the solution of systems of forward-backward stochastic equations of MKV type, though the corresponding forward-backward systems are different. On the one hand, the forward-backward system of the MFG approach appears as a standard forward-backward system derived from the classical form of the stochastic Pontryagin principle, a nonlinear term of MKV type appearing as a result of the fixed point argument of the third step of the MFG approach. On the other hand, as we explained in the special case of scalar interactions derived from \cite{AndersonDjehiche}, the 
forward-backward system associated with the MKV control problem involves additional terms coming from the need to optimize with respect to the interaction terms. 

\vskip 2pt
In this section, we address the solvability of general stochastic forward-backward equations of MKV type. General solvability results are given and then specialized when the coefficients of the forward-backward systems derive either from a MFG or a MKV control problem. Results are stated in a pedagogical manner, almost as a reader's guide to MFG and MKV control problems. In particular, they cover,  in a more abstract fashion, 
the results and computations performed in the linear-quadratic setting (see Sections 
\ref{se:LQMFG}, \ref{se:LQMKV}) and for the other examples discussed in Section \ref{se:examples}. Due to the heavy technical nature of the results and their proofs, we only sketch proofs, all the arguments  being detailed in the forthcoming works \cite{CarmonaDelarue2, CarmonaDelarue3, CarmonaDelarue4} by the two first named authors. While statements are given in dimension 1 only, higher-dimensional versions are available in \cite{CarmonaDelarue2, CarmonaDelarue3, CarmonaDelarue4}.

\subsection{Solvability of General MKV Forward-Backward Stochastic Equations}
\label{subsec:generalMKV}
Forward-backward systems such as \eqref{eq:16:12:1} and \eqref{fo:adjointmkv:0} are here understood as special cases of more general fully coupled forward and backward stochastic differential equations involving the marginal distributions of the forward and backward solutions. Changing the notation ever so slightly in order to accommodate the FBSDEs appearing in both analyses, we consider equations of the form:
\begin{equation}
\label{eq:4:11:1}
\begin{split}
&dX_{t} = B\bigl(t,X_{t},Y_{t},\cL(X_{t},Y_{t})\bigr) dt + \sigma dW_{t}
\\
&dY_{t} = -F\bigl(t,X_{t},Y_{t},\cL(X_{t},Y_{t})\bigr) dt + Z_{t} dW_{t}, \quad 0 \leq t \leq T,
\end{split}
\end{equation}
with $Y_{T} = G(X_{T},\cL(X_{T}))$ as terminal condition. As already mentioned, $\u X$ and $\u Y$ are both one dimensional.

Because of the coupled structure of the system, solvability is a hard question to tackle: fully coupled forward-backward systems are instances of stochastic two-point-boundary-value problems for which both existence and uniqueness are known to fail under standard Cauchy-Lipschitz conditions. When the coefficients do not depend on the marginal distributions of the solutions, the forward and backward equations may be decoupled by taking advantage of the noise. Indeed, when the noise is non-degenerate, it has a decoupling effect by regularizing the underlying FBSDE value function. We refer to \cite{Delarue02} for a review of this strategy. Recall that by FBSDE value function, we mean the function $u$ giving $Y_{t}$ as a function of $X_t$, say $u(t,X_{t})$. When the coefficients $B(t,x,y)$, $F(t,x,y)$ and $G(x)$ are independent of the measure argument, this value function satisfies the quasilinear PDE 
\begin{equation}
\label{eq:16:12:3}
\partial_{t} u(t,x) + \frac{1}{2} \sigma^2 \partial^2_{xx} u(t,x) + B(t,x,u(t,x)) \partial_{x} u(t,x) + 
F(t,x,u(t,x))=0, 
\end{equation}
for $t \in [0,T]$ and $x \in \RR$,
with $u(T,x)=G(x)$. 

When the solution interacts with itself, e.g. with its own distribution, the standard Markov structure breaks down. Indeed, the marginal distribution of the process $\u X$ at time $t$ is needed to compute the transitions towards the future values of the state of the system. Basically, the Markov property must be considered in a larger space, namely the Cartesian product of the state space of the forward process, which is $\RR$ in the present situation, with the space of probability measures on the state space, which is infinite dimensional if the state space is not finite. Put it differently, the relationship between $\u Y$ and $\u X$ is expected to be of the form:
\begin{equation}
\label{eq:1:7:1}
Y_{t} = v\bigl(t,X_{t},{\mathcal L}(X_{t})\bigr), \quad t \in [0,T],
\end{equation}
for some mapping $v : [0,T] \times \RR \times {\mathcal P}_{1}(\RR) \hookrightarrow \RR$. The PDE approach is then not sufficient anymore: the underlying heat kernel derived from the noise still has a smoothing property in the finite-dimensional component, but does not have any regularizing effect in the infinite-dimensional direction of $u$. The regularity of $v$ in $(t,x)$ can be tackled 
by the very simple observation: once the law of $(\u X,\u Y)$ in \eqref{eq:4:11:1} has been computed, 
\eqref{eq:4:11:1} reads as a Markovian FBSDE with $u(t,x) = v(t,x,{\mathcal L}(X_{t}))$ as value function. 

In the forthcoming paper \cite{CarmonaDelarue3}, the solvability of the equation is tackled by a compactness argument and the Schauder fixed point theorem. Because of the expected relationship \eqref{eq:1:7:1} between $\u Y$ and $\u X$, the conditional law of $Y_{t}$ given $X_{t}$ is expected to be a Dirac measure for any $t \in [0,T]$. In other words, the law of the pair $(X_{t},Y_{t})$ is expected to have the form $\varphi(t,\cdot) \diamond \mu_{t}$, where $\varphi(t,\cdot)$ is a continuous mapping from 
$\RR$ into $\RR$, $\mu_{t}$ is a probability measure on $\RR$ and $\varphi(t,\cdot) \diamond \mu_{t}$ is the measure on $\RR^2$, obtained as the image (push--forward) of the distribution $\mu_{t}$ on $\RR$ by the mapping $\RR \ni x \hookrightarrow (x,\varphi(t,x))$. Given 
an element $\mu$ of the space $E={\mathcal P}_{1}(C([0,T]))$ of probability measures  on the space $C([0,T])$ of real valued continuous functions on $[0,T]$ and an element $\varphi$ of the space $F= C_{b}([0,T] \times \RR)$ of real valued bounded continuous functions on $[0,T] \times \RR$, we consider the system
\begin{equation}
\label{eq:4:11:1:b}
\begin{split}
&dX_{t} = B\bigl(t,X_{t},Y_{t}, \varphi(t,\cdot) \diamond \mu_{t} \bigr) dt + \sigma dW_{t}
\\
&dY_{t} = -F\bigl(t,X_{t},Y_{t}, \varphi(t,\cdot) \diamond \mu_{t}\bigr) dt + Z_{t} dW_{t}, \quad 0 \leq t \leq T,
\end{split}
\end{equation}
with $Y_{T} = G(X_{T},\mu_{T})$ as terminal condition. Here $\u\mu=(\mu_t)_{0\le t\le T}$ is the flow of marginal distributions where $\mu_t$ denotes the image (push--forward) of $\mu$ under the $t$-th coordinate map $C([0,T])\ni w\hookrightarrow w(t)\in\RR$.
If this system admits a unique solution, we can consider  as output the measure $\cL(\underline{X})$, which is an element of $E$, 
together with the FBSDE value function linking $\u Y$ with $\u X$, which is an element of $F$ under suitable conditions. The following result is proven in \cite{CarmonaDelarue3}:

\begin{theorem} 
\label{prop:1:12:1}
Assume $\sigma>0$ and that $B$, $F$ and $G$ are bounded Lipschitz continuous with respect to  the space variable for the Euclidean distance, and with respect to the McKean--Vlasov component for the Wasserstein metric, uniformly in the other variables. Then, the mapping $\Phi: (E,F)  \ni (\u\mu,\varphi)  \mapsto (\cL(\underline{X}),u) \in E \times F$, where $u$ is the FBSDE value function such that 
$Y_{t}=u(t,X_{t})$ for any $t \in [0,T]$, has a fixed point.
\end{theorem}
Recall that the Wasserstein metric between two probability measures $\eta$ and $\eta'$ on $\RR^d$, $d \geq 1$, is the square root of the infimum of $\int_{\RR^{2d}} \vert z-z'\vert^2 d \pi(z,z')$ over all the probability measures $\pi$ on $\RR^{2d}$ admitting 
$\eta$ and $\eta'$ as marginals. 

Notice that given any $(\u \mu,\varphi)$ in $E \times F$ , the forward-backward system is standard, so that unique solvability holds, see e.g.  Delarue \cite{Delarue02}. 
The proof consists in showing that $\Phi$ leaves a bounded closed convex subset $\Gamma \subset E \times F$ stable and that the restriction of $\Phi$ to $\Gamma$ is continuous and has a relatively compact range, $E \times F$ being endowed with 
the product of the topology of weak convergence of measures on $E$ and the topology of uniform convergence on compact sets on $F$.
Boundedness of the coefficients plays a crucial role to prove that the probability measures in the range of $\Phi$ are tight.
Positivity of $\sigma$ ensures that the FBSDE value function in the range of $\Phi$ are H\"older continuous and thus live in a compact subset of $F$. A complete proof is given in \cite{CarmonaDelarue3}.  

\vspace{6pt}

The application of Theorem \ref{prop:1:12:1} to the solvability of the adjoint equations in \eqref{eq:16:12:2} and 
\eqref{fo:adjointmkv:0}
 is not so straighforward. Indeed, the coefficients $B$ and $F$ therein derive from a convex Hamiltonian structure, which means that they appear as derivatives with respect to  the state space parameter (and also possibly with respect to  the measure parameter) of the underlying Hamiltonian function. In practice, the Hamiltonian function is expected to grow at infinity at a rate which could be quadratic with respect to the state space parameter and  the control variable, so that boundedness of the coefficients $B$ and $F$ is expected to fail. Similarly, the terminal cost $g$ in \eqref{eq:16:12:2} and  \eqref{fo:adjointmkv:0} may have  quadratic growth, so that the terminal condition in the adjoint equations is not bounded in practice. As a consequence, refinements of Theorem \ref{prop:1:12:1} are necessary for our specific purposes.   

Here is a first straightforward refined statement of this kind: when the coefficient $B$ grows at most linearly in  $(x,y)$ and $\cL(X_{t},Y_{t})$, and $F$ and $G$ are bounded, existence of a solution still holds true. When we say that $B$ is at most of linear growth with respect to $(x,y,\cL(X_{t},Y_{t}))$, we mean that
\begin{equation*}
|B(t,x,y,\cL(X_{t},Y_{t}))| \leq C \biggl[ 1 + |x| + |y| + \biggl( \int_{\RR^{2}}
|x'|^2 d \mu(x') \biggr)^{1/2} \biggr], \quad \mu = \cL(X_{t},X_{t}).
\end{equation*}
Solvability follows from a compactness argument again. Approximating the drift 
$B$ by a sequence of  bounded drifts $(B_n)_{n \geq 1}$, we can prove that, for the corresponding solutions $((X^n_t,Y^n_t=u^n(t,X^n_t))_{0 \leq t \leq T})_{n \geq 1}$, the distributions $(\cL(\underline{X}^n))_{n \geq 1}$ are tight and that the functions $(u^n)_{n \geq 1}$ are equicontinuous on compact subsets of 
$[0,T] \times \RR$. The proof follows from two key observations: first, since $f$ and $g$ are bounded, the functions $(u^n)_{n \geq 1}$ are uniformly bounded, so that 
the processes $(\underline{X}^n)_{n \geq 1}$ are tight; second, by smoothing properties of keat kernels, the functions $(u^n)_{n \geq 1}$ are locally uniformly continuous, see \cite{Delarue03}. Extracting a convergent subsequence, we can pass to the limit by using stability of forward-backward stochastic differential equations. 

Here is another straightforward refinement: when $G$ is bounded, and $F$ is bounded with respect to all the parameters but $y$, and has a linear growth with respect to  $y$, the solvability still holds. Indeed, a standard maximum principle for BSDEs says that 
the process $\u Y$ is bounded, the bound depending upon the bound of $G$ and the growth of $F$ only, so that $F$ may be seen as a bounded driver.

\subsection{Solvability of the MFG Adjoint Equations}
In the MFG framework, the FBSDE system is given by \eqref{eq:16:12:2}. With the notation \eqref{eq:4:11:1}, it reads
\begin{equation*}
\begin{split}
&B\bigl(t,X_{t},Y_{t},{\mathcal L}(X_{t},Y_{t})\bigr)
= b \bigl(t,X_{t},{\mathcal L}(X_{t}),\hat{\alpha}^{{\mathcal L}(X_{t})}(t,X_{t},Y_{t}) \bigr),
\\
&F\bigl(t,X_{t},Y_{t},{\mathcal L}(X_{t},Y_{t})\bigr)
=  \partial_{x} H^{{\mathcal L}(X_{t})} \bigl(t,X_{t},Y_{t},\hat{\alpha}^{{\mathcal L}(X_{t})}(t,X_{t},Y_{t}) \bigr),
\\
&G\bigl(X_{T},{\mathcal L}(X_{T})\bigr) = \partial_{x} g\bigl(X_{T},{\mathcal L}(X_{T})\bigr),
\end{split}
\end{equation*}
where $H$ and $\hat{\alpha}$ are given by \eqref{fo:hmu} and \eqref{fo:alphahat}. We emphasize that the dependence upon the McKean--Vlasov interaction is limited to the law of $X_{t}$ only and that the coefficients do not involve the law of $Y_{t}$. When 
$\partial_{x}g$ is bounded and $\partial_{x} H$ is bounded with respect to all the parameters but $y$, and is at most of linear growth  in $y$, Theorem \ref{prop:1:12:1} applies directly, provided that the coefficients are Lipschitz continuous with respect to $x$, $y$ and 
$\mu$, the Lipschitz continuity with respect to $\mu$ being with respect to  the Wasserstein distance. 
 (See \cite[Section 3]{Cardaliaguet} for the analytical counterpart of this result.)
In particular,  
the minimizer $\hat{\alpha}^{{\mu_{t}}}(t,x,y)$ must be Lipschitz continuous as well. In the examples tackled below, the Hamiltonian is always strictly convex in the variable $\alpha$ so that the regularity of the minimizer follows from the implicit function theorem. 

In practice, the cost functions $g$ and $H$ are expected to be of quadratic growth in $x$ so that $\partial_{x} g$ and 
$\partial_{x} H$ cannot be bounded. The main point is thus to allow $F$ and $G$ to grow at most linearly in the statement of 
Theorem \ref{prop:1:12:1}. The strategy developed in \cite{CarmonaDelarue2} consists in approximating the 
cost functions of the MFG control problem by a sequence of cost functions with bounded derivatives in $x$ preserving the convexity of the Hamiltonian in $(x,\alpha)$. Here is the solvability result obtained in \cite{CarmonaDelarue2}:

\begin{theorem}
\label{thm:adjoint:1}
Let us assume that

(i) The cost function $f$ is Lipschitz continuous in the variables $x$, $\mu$ and $\alpha$ on any balls of center $0$ and radius 
$R$, that is on sets where $x$ and $\alpha$ are bounded by $R$ and the second-order moment of $\mu$ is by bounded by $R^2$, the Lipschitz constant being at most of linear growth in $R$. Moreover, $f$ is twice-continuously differentiable with respect to 
$x$ and $\alpha$, with uniformly bounded second-order derivatives (uniformly in $(t,x,\mu,\alpha)$). The partial derivative  $\partial_{\alpha} f(t,x,\mu,0)$ is uniformly bounded and the product $
x \partial_{x} f(t,0,\delta_{x},0)$ is always non-negative. (Here, $\delta_{x}$ denotes the Dirac mass at point $x$.)
Moreover, the function $f$ is convex in  $(x,\alpha)$ for $(t,\mu)$ fixed in the sense that there exists a constant $\lambda >0$ such that:
\begin{equation*}
f(t,x',\mu,\alpha') - f(t,x,\mu,\alpha) - \langle (x'-x,\alpha'-\alpha),\nabla_{(x,\alpha)} f(x,\mu,\alpha) \rangle \geq \lambda 
\vert \alpha' - \alpha \vert^2;
\end{equation*} 

(ii) The terminal cost $g$ is Lipschitz continuous in the variables $x$ and $\mu$ on any balls of center $0$ and radius 
$R$ (in the same sense as $f$)  and the Lipschitz constant is at most of linear growth in $R$. Moreover, $g$ 
is twice-continuously differentiable with respect to  $x$, with uniformly bounded second-order derivatives 
(uniformly in $(x,\mu)$), $g$ is convex in $x$ when $\mu$ is fixed and the product $x \partial_{x} g(0,\delta_{x})$ is always non-negative;
 
(iii) $b$ is affine in $x$ and $\alpha$ in the sense that $b(t,x,\mu,\alpha) = b_{0}(t,\mu)+b_{1}(t)x+b_{2}(t) \alpha$ 
where $b_{0}$, $b_{1}$ and $b_{2}$ are bounded and $b_{0}$ is Lipschitz-continuous with respect to  
$\mu$ for the Wasserstein distance.

Then, the adjoint equation \eqref{eq:16:12:2} has at least one solution.
\end{theorem}

There are three types of assumptions in the statement of Theorem \ref{thm:adjoint:1}. First, by linearity of $b$ in $(x,\alpha)$ and by convexity of $f$ in $(x,\alpha)$, the Hamiltonian $H$ is convex in $(x,\alpha)$. It is even strictly convex in $\alpha$. 
Following \eqref{eq:16:12:2}, this says that the stochastic Pontryagin principle applies in the MFG procedure. Second, Lipschitz properties say that the coefficients are regular with respect to the parameters $x,\alpha$ and $\mu$. Actually, by a careful inspection of $(i)$ and $(ii)$, the reader might notice that nothing is said about the regularity of $\partial_{x} H$ and 
$\partial_{x} g$ in the parameter $\mu$, whereas $F$ and $G$ are asked to be Lipschitz continuous in $\mu$ in the statement 
of Theorem \ref{prop:1:12:1}. Indeed, as proven in \cite{CarmonaDelarue2}, only the regularity of the cost functions with respect to the parameter $\mu$ really matters within the framework of control theory. Third, the sign conditions on 
$x \partial_{x} f(t,0,\delta_{x},0)$ and $
x \partial_{x} g(t,0,\delta_{x},0)$ sound as mean-reverting assumptions of a weak type: when approximating the coefficients by bounded ones, the sign conditions guarantee that the expectations of the approximating solutions do not blow up. 

In comparison with Section 
\ref{se:LQMFG}, we emphasize that the mean-field linear-quadratic models investigated therein do not satisfy the boundedness condition required for $b_{0}$ in $(iii)$. In the current setting, the boundedness condition ensures that the fixed point measures of the approximating MFG problems are tight. In the linear-quadratic framework, the tightness condition can be expressed differently, as the problem may be reformulated directly, by investigating the solvability of associated Riccati equations. Moreover, the current sign conditions must be compared with the sign conditions in the statement of Theorem 
\ref{thm:30:12:1}. In the framework of Theorem 
\ref{thm:30:12:1}, $x \partial_{x} g(0,\delta_{x}) = q \o q x^2$ and 
$x \partial_{x} f(t,0,\delta_{x},0) = m_{t} \o m_{t} x^2$, so that the sign condition in Theorem \ref{thm:adjoint:1} hold if $q,\o q \geq 0$
and $\u m, \u {\o m} \geq 0$ respectively. Keep in mind that we can always assume $q \geq 0$ and 
$\u m \geq 0$. Clearly, this is stronger than the sign conditions $q(q+\o q) \geq 0$ and $\inf_{t \in [0,T]} [m_{t}(m_{t}+\o m_{t})] \geq 0$ required in Theorem 
\ref{thm:30:12:1}. 
\vspace{5pt}

\textit{Sketch of Proof (Theorem \ref{thm:adjoint:1}).} As already announced, the proof consists in 
approximating $f$ and $g$ by two sequences $(f^n)_{n \geq 1}$ and $(g^n)_{n \geq 1}$ satisfying the same assumptions, but with bounded derivatives in $x$. The construction of the approximating sequence follows from arguments of convex analysis, which may be found in \cite{CarmonaDelarue2}. We use the convexity of $f$ and $g$ to 
express $f$ and $g$ as their own Legendre bi-conjugates and then truncate the dual variable in the Legendre representation.

By Theorem \ref{prop:1:12:1}, we then expect\footnote{We say ``expect'' only since we said nothing about the smoothness of 
$\partial_{x} H$ and $\partial_{x} g$ with respect to $\mu$. We refer to \cite{CarmonaDelarue2} for the complete argument.} to find, for each $n \geq 1$, a solution $(X^n,Y^n)$ to the adjoint equations associated with the approximate Hamiltonian
\begin{equation*}
H^{n,\mu}(t,x,y,\alpha) = b(t,x,\mu,\alpha) y + f^{n}(t,x,\mu,\alpha), \quad t \in [0,T], \ x,y,\alpha \in \RR, \ \mu \in {\mathcal P}_{1}(\RR),
\end{equation*}
and therefore with the coefficients
$B^n(t,x,y,\mu)=b(t,x,\mu,\hat{\alpha}^n(t,x,y,\mu))$, 
$F^n(t,x,y,\mu)=\partial_{x}f^n(t,x,\mu,\hat{\alpha}^n(t,x,y,\mu)) + b_{1}(t) y$ and
$G^n(x,\mu)=\partial_{x}g^n(x,\mu)$,
where 
\begin{equation*}
\hat{\alpha}^n(t,x,y,\mu) = \textrm{argmin}_{\alpha \in \RR} H^{n,\mu}(t,x,y,\alpha).
\end{equation*}
We also expect $\hat{\alpha}^n(t,x,y,\mu)$ to converge towards $\hat{\alpha}(t,x,y,\mu)$ as $n \rightarrow + \infty$.

For $n \geq 1$, we denote by $\underline{X}^n$ the process controlled by $(\hat{\alpha}^n_{t} = \hat{\alpha}^n(t,X_{t}^n,Y_{t}^n,
\cL(X_{t}^n)))_{0 \leq t \leq T}$ and by $u^n$ the FBSDE value function such that $Y^n_{t}=u^n(t,X^n_{t})$, for any 
$t \in [0,T]$. Our goal is to establish tightness of the processes $(\underline{X}^n)_{n \geq 1}$ and relative compactness of the functions $(u^n)_{n \geq 1}$ for the topology of uniform convergence on compact subsets of $[0,T] \times \RR$. The first step is to prove that 
\begin{equation}
\label{eq:28:12:5}
\sup_{n \geq 1} \EE \biggl[ \int_{0}^T \vert \hat{\alpha}^n_{s} \vert^2 ds \biggr] < + \infty, 
\end{equation}
and the second step is to prove that the growth of the functions $(u^n)_{n \geq 1}$ can be controlled, uniformly in $n \geq 1$. 

In order to prove \eqref{eq:28:12:5}, we take advantage of the Hamiltonian structure of the coefficients by comparing the behavior of 
$\underline{X}^n$ to the behavior of a \textit{reference} controlled process, driven by specific values of the control. Below, the generic notations for the reference controlled process and for the reference control are $\underline{U}^n$ and 
$\underline{\beta}^n$ respectively. Two different versions are considered for $\underline{U}^n$: they are driven by the controls $\underline{\beta}^n \equiv 0$ and $({\beta}^n_{s} = \EE(\hat{\alpha}_{s}^n))_{0 \leq s \leq T}$, respectively. For each of these controls, we compare the corresponding cost to the optimal one by using the stochastic Pontryagin principle and, subsequently, to derive some information on the optimal control $(\hat{\alpha}_{s}^n)_{0 \leq s \leq T}$. The starting point is to compare the behavior of $\underline{X}^n$ to the behavior of the process $\underline{U}^n$ controlled by the deterministic control $(\beta^n_{s} = {\mathbb E}(\hat{\alpha}^n_{s}))_{0 \leq s \leq T}$. As explained right below, the goal is to obtain a uniform bound for the $L^2$-norm in time of the variance of $(\hat{\alpha}^n_{s})_{0 \leq s \leq T}$. Indeed, because of the specific affine structure of $b$ in $(iii)$, it can be checked that 
$\EE(U^n_{s}) = {\mathbb E}(X^n_{s})$, for $0 \leq s \leq T$, and that $\sup_{n \geq 1} \sup_{0 \leq s \leq T} \textrm{Var}(U^n_{s}) < + \infty$. By the stochastic Pontryagin principle and by the growth and convexity properties in 
$(i)$--$(iii)$, we then claim
\begin{equation*}
\begin{split}
&g^n\bigl(\EE(X^n_{T}),\cL(X^n_{T})\bigr) + \int_{t}^T 
\bigl[ \lambda  \textrm{Var} (  \hat{\alpha}_{s}^n ) +
f^n \bigl(s,\EE(X_{s}^n),\cL(X_{s}^n),\EE(\hat{\alpha}_{s}^n) \bigr) \bigr] 
ds 
\\
&\hspace{15pt}
\leq \EE \biggl[ g^n\bigl(U^n_{T},\cL(X^n_{T})\bigr) +  \int_{t}^T f^n \bigl(s,U_{s}^n,\cL(X_{s}^n),0\bigr) 
ds \biggr],
\end{split}
\end{equation*}
from which we can derive (see \cite{CarmonaDelarue2})
\begin{equation}
\label{eq:28:12:6}
\sup_{n \geq 1} \biggl[ \int_{0}^T \textrm{Var}(\hat{\alpha}_{s}^n) ds + 
  \sup_{0 \leq s \leq T} \textrm{Var}(X^n_{s})
\biggr] 
\leq 
c \biggl( 1 + \EE \int_{0}^T \vert \hat{\alpha}_{s}^n \vert^2 ds \biggr)^{1/2}.
\end{equation}
Bound \eqref{eq:28:12:6} says that \eqref{eq:28:12:5} holds  if the expectations of the variables $(X^n_{t})_{0 \leq t \leq T}$ are bounded. In other words, we can only focus on the expectations of the variables $(X^n_{t})_{0 \leq t \leq T}$. 

The second step in the proof of \eqref{eq:28:12:5} consists in comparing the behavior of $\underline{X}^n$ to the behavior of the process controlled by the null control. Since no confusion is possible, we still denote by $\underline{U}^n$ the process controlled by the null control $\underline{\beta}^n \equiv 0$. By boundedness of $b_{0}$ in $(iii)$, it is easily checked that $\sup_{n \geq 1} {\mathbb E}[\sup_{0 \leq s \leq T} \vert U^n_{s} \vert^2] < + \infty$. Applying the stochastic Pontryagin principle again together 
with Jensen's inequality,
\begin{equation}
\label{eq:29:12:15}
\begin{split}
&g^n\bigl(\EE(X^n_{T}),\cL(X^n_{T})\bigr) + \int_{t}^T 
\bigl[ \lambda  \EE \bigl( \vert \hat{\alpha}_{s}^n \vert^2\bigr) +
f^n \bigl(s,\EE(X_{s}^n),\cL(X_{s}^n),\EE(\hat{\alpha}_{s}^n) \bigr) \bigr] 
ds 
\\
&\hspace{15pt}
\leq \EE \biggl[ g^n\bigl(U^n_{T},\cL(X^n_{T})\bigr) +  \int_{t}^T f^n \bigl(s,U_{s}^n,\cL(X_{s}^n),0\bigr) 
ds \biggr].
\end{split}
\end{equation}
The main point in \cite{CarmonaDelarue2} is to deduce that
\begin{equation*}
\EE(X^n_{T}) \partial_{x} g^n\bigl(0,\delta_{\EE(X^n_{T})}\bigr) + \int_{t}^T 
\bigl[ \frac{\lambda}{4}  \EE \bigl( \vert \hat{\alpha}_{s}^n \vert^2\bigr) +
\EE(X^n_{s})
\partial_{x} f^n \bigl(s,0,\delta_{\EE(X^n_{s})},0 \bigr) \bigr] 
ds 
\leq c.
\end{equation*}
 By the sign property in $(i)$ and $(ii)$, \eqref{eq:28:12:5} follows. We deduce that
$\EE[\sup_{0 \leq s \leq T} \vert X^n_{s} \vert^2] < + \infty$. We also deduce that the processes $(\underline{X}^n)_{n \geq 1}$ are tight. 

Once we have a uniform estimate for the moments of the $(\underline{X}^n)_{n \geq 1}$, we are led back to FBSDEs of the standard Markovian type, each of them appearing as the adjoint system of a standard optimal stochastic control problem. Standard FBSDE arguments then imply that 
$\vert u^n(t,x) \vert \leq c(1+ \vert x \vert)$. Uniform continuity of the family $(u^n)_{n \geq 1}$ on compact subsets follows from the smoothing effect of heat kernels, see \cite{Delarue03}. 

Let $(u,\u\mu)$ be the limit of a convergent subsequence which we still denote $(u^n,\cL(\underline{X}^n))_{n \geq 1}$. The smoothing effect of heat kernels can also be used to prove that the functions $(u^n)_{n \geq 1}$ are uniformly locally Lipschitz-continuous in space, so that $u$ itself is locally Lipschitz-continuous in space. Since $u$ grows at most linearly, the SDE
\begin{equation*}
dX_{t} = B\bigl(t,X_{t},u(t,X_{t}),\mu_{t}\bigr) dt + \sigma d
 W_{t},
\end{equation*}
is strongly solvable. We can also solve the backward equation
\begin{equation*}
dY_{t} = -F\bigl(t,X_{t},Y_{t},\mu_{t}\bigr) dt + Z_{t} d W_{t}, \quad Y_{T} = G(X_{T},\mu_{T}).
\end{equation*}
By stability of forward and backward equations separately, it is plain to see that $\underline{X}$ is the limit of the $(\underline{X}^n)_{n \geq 1}$ and that 
$\underline{Y}$ is the limit of the $((u^n(t,X^n_{t}))_{0 \leq t \leq T})_{n \geq 1}$. Passing to the limit in the approximating forward-backward systems, we conclude that $(\underline{X},\underline{Y})$ solves the adjoint equation in \eqref{eq:16:12:2}. \qed

\subsection{Solvability of the adjoint Equations of the MKV Control Problem}
In the framework of the optimal control of MKV processes, the forward-backward system is given by 
\eqref{fo:adjointmkv:0}. As for the MFG approach, it can be recast as a system of the type 
\eqref{eq:4:11:1}, the coefficients of the backward component depending on the law of $Y$ through the term $\EE\{\partial_{x'}b(t,X_t,\EE\{\psi(X_t)\},\alpha_t)Y_t\}\partial_x\psi(X_t)$. We claim:

\begin{theorem}
\label{th:sol:mkv}
In the framework of Theorem \ref{th:pontryaginmkv}, assume that $\gamma$ and $\zeta$ are twice-differentiable convex functions, with bounded second-order derivatives and that $b$, $f$ and $g$ are twice-differen- -tiable functions with bounded second-order derivatives that satisfy  (A1) - (A3).  Assume in addition that $b$ is linear with respect to $x$, $x'$ and $\alpha$, namely 
$b(t,x,x',\alpha)=b_{0}(t) + b_{1}(t) x + b_{2}(t) x' + b_{3}(t) \alpha$, where $b_{0}$, $b_{1}$, $b_{2}$ and $b_{3}$ are bounded functions. Assume furthermore that $\psi(\xi)=\xi$ so that $\langle \psi,\mu \rangle$ stands for the expectation of $\mu$. Assume finally that 
$f$ is $\lambda$-convex with respect to  $\alpha$, for some $\lambda >0$, so that $H$ in \eqref{fo:hamiltonianmkv} is also 
$\lambda$-convex with respect to  $\alpha$. 

Then, the forward-backward system \eqref{fo:adjointmkv:0}, with $\alpha_{t} = \hat{\alpha}(t,X_{t},Y_{t},{\mathcal L}(X_{t}))$, is uniquely solvable, where $\hat{\alpha}(t,x,y,\mu) = \argmin_{\alpha} H(t,x,\langle \psi,\mu\rangle,\langle \gamma,\mu \rangle,y,\alpha)$.
\end{theorem}

\textit{Sketch of Proof.} We first notice that the function $\hat{\alpha}$ is well-defined and Lipschitz-continuous w.r.t to  $x$, $y$ and $\mu$. As in the proof of Theorem \ref{thm:adjoint:1}, we consider sequences of convex functions $(f^n)_{n \geq 1}$, $(g^n)_{n \geq 1}$, $(\gamma^n)_{n \geq 1}$ and $(\zeta^n)_{n \geq 1}$, with bounded first and second-order derivatives, that converge towards $f$, $g$, $\gamma$ and $\zeta$ respectively. As explained in \cite{CarmonaDelarue4}, we can choose 
$(f^n)_{n \geq 1}$, $(g^n)_{n \geq 1}$, $(\gamma^n)_{n \geq 1}$ and $(\zeta^n)_{n \geq 1}$ to be non-decreasing sequences.
For each $n \geq 1$, we define the Hamiltonian $H^n$ 
by replacing $f$ by $f^n$ in \eqref{fo:hamiltonianmkv}, and we consider the forward-backward system
\eqref{fo:adjointmkv:0} by replacing $(f,g,\gamma,\zeta)$ therein by $(f^{n},g^{n},\gamma^{n},\zeta^{n})$
and $\alpha_{t}$ by  $\hat{\alpha}^n_{t}=\hat{\alpha}^n(t,X_{t},Y_{t},{\mathcal L}(X_{t}))$, the argument of the minimum of $H^n(t,x,\langle \psi,\mu\rangle,\langle \gamma,\mu \rangle,y,\alpha)$: it is the forward-backward system associated with the Hamiltonian $H^n$. We denote by $(\u X^n,\u Y^n)$ its solution.

By the stochastic Pontryagin principle, the process $(\hat{\alpha}^n_{t})_{0 \leq t \leq T}$ minimizes the cost function 
\eqref{fo:mkvobjective}, with $f$, $g$, $\gamma$ and $\zeta$ therein replaced by $f^n$, $g^n$, $\gamma^n$ and $\zeta^n$. In particular, $J^n(\underline{\hat{\alpha}}^n) \leq J^{n}(\underline{\hat{\alpha}}^{n+p})$, for $p \geq 0$. 
Actually, by $\lambda$-convexity of $f^n$ with respect to  $\alpha$, we can even write, as a variant of the stochastic Pontryagin principle,
\begin{equation*}
J^n\bigl(\underline{\hat{\alpha}}^n\bigr) + \lambda {\mathbb E} \int_{0}^T 
 \vert \hat{\alpha}_{s}^n - \hat{\alpha}_{s}^{n+p} \vert^2 ds 
\leq J^n\bigl(\underline{\hat{\alpha}}^{n+p}\bigr).
\end{equation*}
Since the sequences $(f^n)_{n \geq 1}$ and $(g^n)_{n \geq 1}$ are non-decreasing, 
we deduce that, for any admissible control $\underline{\beta}$, $J^n(\underline{\beta}) \leq J^{n+p}(\underline{\beta})$. Therefore,
\begin{equation*}
J^n\bigl(\underline{\hat{\alpha}}^n\bigr) + \lambda {\mathbb E} \int_{0}^T 
\vert \hat{\alpha}_{s}^n - \hat{\alpha}_{s}^{n+p} \vert^2 
ds 
\leq J^{n+p}\bigl(\underline{\hat{\alpha}}^{n+p}\bigr).
\end{equation*}
In particular, the sequence $(J^n(\underline{\hat{\alpha}}^n))_{n \geq 1}$ is non-decreasing. Thus, it has a limit since it is bounded by 
$\sup_{n \geq 1}(J^n(0))<\infty$. We deduce that 
$(\underline{\hat{\alpha}}^n)_{n \geq 1}$ is a Cauchy sequence. By stability of standard McKean Vlasov SDEs, we deduce that  that
\begin{equation*}
\lim_{n \rightarrow + \infty}
\sup_{p \geq 1}{\mathbb E} \bigl[ \sup_{0 \leq s \leq T} \vert X_{s}^{n+p} - X_{s}^n \vert^2 \bigr] = 0.
\end{equation*}
This proves that the processes $(\underline{X}^n)_{n \geq 1}$ converge for the norm ${\mathbb E}[\sup_{0 \leq s \leq T} | \cdot |^2 ]^{1/2}$, towards a continuous adapted process
$\underline{X}$. By standard results of stability for BSDEs, we then deduce that the same holds for the processes $(\underline{Y}^n)_{n \geq 1}$, that is
\begin{equation*}
\lim_{n \rightarrow + \infty}
\sup_{p \geq 1}{\mathbb E} \bigl[ \sup_{0 \leq s \leq T} \vert Y_{s}^{n+p} - Y_{s}^n \vert^2 \bigr] = 0.
\end{equation*}
One then deduce the existence of the limit $\underline{Y}$, and the fact that $(\underline{X},\underline{Y})$ satisfies the forward-backward system \eqref{fo:adjointmkv:0}. 
Uniqueness easily follows from the stochastic Pontryagin principle. \qed

\begin{remark}
Uniqueness in Theorem 
\ref{th:sol:mkv} says that, for any $t \in [0,T]$, $Y_{t}$ is a function of $X_{t}$ and of the law of $X_{t}$, as in \eqref{eq:1:7:1}. This also follows from a standard change of filtration, as usually done in the theory of BSDEs. In particular, the optimal control $\underline{\hat{\alpha}}$ is a feedback control, as expected.  
\end{remark}

\bibliographystyle{plain}

\end{document}